\documentclass[11pt,twoside]{article}
\usepackage{url}
\usepackage{courier}
\usepackage{manfnt}
\usepackage{relsize}
\usepackage{setspace}
\usepackage{stmaryrd}
\usepackage{amsmath}
\usepackage{latexsym}
\usepackage{psfrag}
\usepackage{amsfonts} 
\usepackage{amsmath} 
\usepackage{amssymb}
\usepackage{bbm}
\usepackage{amsthm} 
\usepackage{mathrsfs} 
\usepackage{comment}
\usepackage{enumerate}
\usepackage{url}
\usepackage{graphicx}

\usepackage{rotating}
\usepackage[config, labelfont={normalsize,bf}, textfont=normalsize]{caption,subfig}
\usepackage{epsfig}
\usepackage[stable]{footmisc}
\usepackage{endnotes} 
\usepackage[english]{babel}
\usepackage{makecell}
\usepackage{rotating}
\usepackage{multirow}
\usepackage[hmargin=1in,vmargin=1in]{geometry}
\usepackage{mathtools}
\theoremstyle{definition}

\newtheorem{shortassumption}{A\!\!}

\theoremstyle{plain}
\newtheorem{theorem}{Theorem}

\newtheorem{corollary}{Corollary}
\newtheorem{lemma}{Lemma}

\theoremstyle{definition}

\newcommand{\R}{\mathbb{R}}

\renewcommand{\hat}[1]{\widehat{#1}}

\newcommand{\I}[1]{\mathbbm{1}}

\newcommand{\D}{\mathcal{D}}

\newcommand{\F}{\mathbb{F}}

\newcommand{\X}{\mathcal{X}}

\newcommand{\tsum}{\textstyle\sum}

\renewcommand{\P}{\mathbb{P}}
\newcommand{\E}{\mathbb{E}}

\newcommand{\ve}{\varepsilon}

\newcommand{\var}{\operatorname{var}}

\newcommand{\ts}{\textstyle}

\newcommand{\err}{\text{err}}
\newcommand{\Err}{\textsc{Err}}

\newcommand{\mnorm}[1]{\left\vert\kern-1.5pt\left\vert\kern-1.5pt\left\vert #1\right\vert\kern-1.5pt\right\vert\kern-1.5pt\right\vert}

\newcommand\wcline[1]{%
    \noalign{\xdef\origarrayrulewidth{\the\arrayrulewidth}%
    \global\arrayrulewidth 5\arrayrulewidth}%
    \cline{#1}%
    \noalign{\global\arrayrulewidth\origarrayrulewidth}%
}

\begin{document}
\begin{center}

{\bf{\LARGE{Estimating a sharp convergence bound for randomized ensembles}}}

%

\vspace{.2cm}
{\large{
\begin{center}
Miles E. Lopes\footnote{ Department of Statistics, University of California, One Shields Avenue, Davis, CA 95616, USA, email: {\tt{melopes@ucdavis.edu}}. This research was partially supported by NSF grant DMS 1613218. The author thanks Ethan Anderes, James Sharpnack, and Philip Kegelmeyer for helpful feedback.} \\[1ex]
\date{}
\end{center}
}}
\vspace{.1in}


\begin{abstract}
When randomized ensembles such as bagging or random forests are used for binary classification, the prediction error of the ensemble tends to decrease and stabilize as the number of classifiers increases. However, the precise relationship between prediction error and ensemble size is unknown in practice. In the standard case when classifiers are aggregated by majority vote, the present work offers a way to quantify this convergence in terms of \mbox{``algorithmic variance,''} i.e.~the variance of prediction error due only to the randomized training algorithm. Specifically, we study a theoretical upper bound on this variance, and show that it is sharp --- in the sense that it is attained by a specific family of randomized classifiers. Next, we address the problem of estimating the unknown value of the bound, which leads to a unique twist on the classical problem of non-parametric density estimation. In particular, we develop an estimator for the bound and show that its MSE matches optimal non-parametric rates under certain conditions. (Concurrent with this work, some closely related results have also been considered in Cannings and Samworth~\cite{cannings2017} and Lopes~\cite{lopes2019}.) 
\end{abstract}

\end{center}

\section{Introduction}
During the past two decades, randomized ensemble methods such as bagging and random forests have become established as some of the most popular prediction methods~\cite{breiman1996, breiman2001}. Although the literature has thoroughly explored how the prediction error of these methods depends on training sample size (e.g.~\cite{buhlmannyu,breimanconsistency,hallsamworth,linjeon,biau2008,biau2012,scornetconsistency} among others),
comparatively little is known about how the prediction error depends on the number of classifiers (ensemble size). In particular, only a handful of works have considered this question from a theoretical standpoint~\cite{NgJordan,hernandez2013,cannings2017,lopes2019} (cf.~Section~\ref{sec:related}).

Indeed, it is of basic interest to have some guarantee that an ensemble is large enough so that it has reached ``algorithmic convergence''  --- i.e.~when the prediction error is close to the ideal level of an infinite ensemble (on a given dataset). 
Specifically, this type of guarantee prevents wasteful computation on an excessively large ensemble, and it ensures that any potential gains in accuracy from a larger ensemble are minor.

\subsection{Background and setup}
At a high level, ensemble methods are applied to binary classification in the following way. Given a set of $n$ labeled training samples $\D:=\{(X_j,Y_j)\}_{j=1}^{n}$ in a generic sample space $\X\times \{0,1\}$, an algorithm is used to train an ensemble of base classifiers \mbox{$Q_i :\X\to \{0,1\}$}, $i=1,\dots, t$. The predictions of the classifiers are then aggregated by a particular rule, with majority vote being the standard choice for bagging and random forests.
More precisely, if a random point $(X,Y)$ is sampled from $\X\times \{0,1\}$, with the label $Y$ being unknown, then we  write the majority vote of the ensemble as a binary indicator function \mbox{$M_t(X):=1\{\bar{Q}_t(X)\geq \frac{1}{2}\}$,} where $\bar{Q}_t(\cdot):=\frac{1}{t}\sum_{i=1}^t Q_i(\cdot)$. Also, for simplicity, we assume going forward that $t$ is odd to eliminate the possibility of ties.

\paragraph{Randomized ensembles} Since bagging and random forests are the motivating examples for our analysis, we now review how randomization arises in these methods. In the case of bagging, the method generates a collection of random sets $\D_1^*, \dots,\D_t^*$, each of size $n$, by sampling with replacement from $\D$. In turn, for each $i=1,\dots,t$, a classifier $Q_i$ is trained on $\D_i^*$, using a ``base classification method'' (such as a decision tree method). Similarly, random forests may be viewed as an extension of bagging, since it generates the sets $\{\D_i^*\}$ in the manner above, but adds one extra ingredient: a randomized feature selection rule for training each $Q_i$ on $\D_i^*$. (We refer to the book~\cite[Ch. 15]{elements} for further details.)

An important commonality of bagging and random forests is that the classifiers $\{Q_i\}$  can be represented in the following way. Specifically, there is a deterministic function, say $g$, and a sequence of i.i.d.~``randomizing parameters'' $\xi_1,\xi_2,\dots$, independent of $\D$, such that each classifier $Q_i$ can be written as
\begin{equation}\label{rfrep}
Q_i(x)=g(x,\D,\xi_i), \ \ \  \ \ \text{ for all } x\in\mathcal{X}.
\end{equation}
In addition to making the presence of algorithmic randomness explicit through the variables $\{\xi_i\}$, this representation also makes it clear that the classifiers $\{Q_i\}$ have the property of being conditionally i.i.d., given $\D$.
In the case of bagging, each $\xi_i$ plays the role of the random set $\D_i^*$, whereas in the case of random forests, each $\xi_i$ encodes both the set $\D_i^*$ as well as randomly selected sets of features~\cite[cf.~Definition~1.1]{breiman2001}. So, as a way of unifying our results, we will consider a general class of ensembles of this form, per the following assumption.
\begin{shortassumption}\label{assump:rep}
The sequence of randomized classifiers $\{Q_i\}$ can be represented in the form~\eqref{rfrep}.
\end{shortassumption}

In addition to bagging and random forests, this assumption is satisfied by the voting Gibbs classifier~\cite{NgJordan}, as well as ensembles generated using random projections, e.g.~\cite{cannings2017}. In the case of the voting Gibbs classifier, a sequence of $t$ i.i.d.~binary samples is drawn from a posterior distribution, and then predictions are made with the majority vote of these samples. Alternatively, for methods based on random projections, a sequence of $t$ i.i.d.~random matrices are used to create $t$ projected versions of $\D$, and in turn, these different versions of $\D$ are used to train an ensemble of classifiers that are aggregated by majority vote. Furthermore, there are several other relatives of bagging and random forests that satisfy\textbf{A1}, including those in~\cite{tinkamho1998,dietterich2000,buhlmannyu}.
Lastly, it is worth pointing out that boosting methods typically do not satisfy \textbf{A1}, and we refer to the book~\cite{schapire2012boosting} for an overview of algorithmic convergence in that context.

\paragraph{Error rates} \ \
 We now define the error rates that will be the focus of our analysis. Letting $l\in\{0,1\}$ be a placeholder for the class label, we define $\mu_l:=\mathcal{L}(X|Y=l)$ as the  distribution of the test point $X$, given that it is drawn from class $l$. Also let $\boldsymbol \xi_t=(\xi_1,\dots,\xi_t)$.
Then, for a particular realization of the classifiers $\{Q_i\}$, trained on a particular set $\D$, the class-wise prediction error rates are defined by
\begin{equation}\label{errdef}
\small
\begin{split}
\Err_{t,l} &:= \int_{\mathcal{X}} 1\big\{|\bar{Q}_t(x)-l|\geq \ts\frac{1}{2}\big\} d\mu_l(x) \ = \ \P\Big(M_t(X)=1-l \, \Big| \, \boldsymbol \xi_t, \D, Y=l\Big).
\end{split}
\end{equation}
In this definition, it is important to emphasize that $\Err_{t,l}$ is a random variable. Specifically, there are two sources of randomness in $\Err_{t,l}$, which are the algorithmic randomness from the variables $\boldsymbol \xi_t$, and the randomness from the data $\D$.
 Since our goal is to quantify algorithmic convergence on a given set $\D$, we will always analyze $\Err_{t,l}$ conditionally on $\D$.

\begin{figure}[h!]
\centering
{\includegraphics[angle=0,
  width=0.45\linewidth]{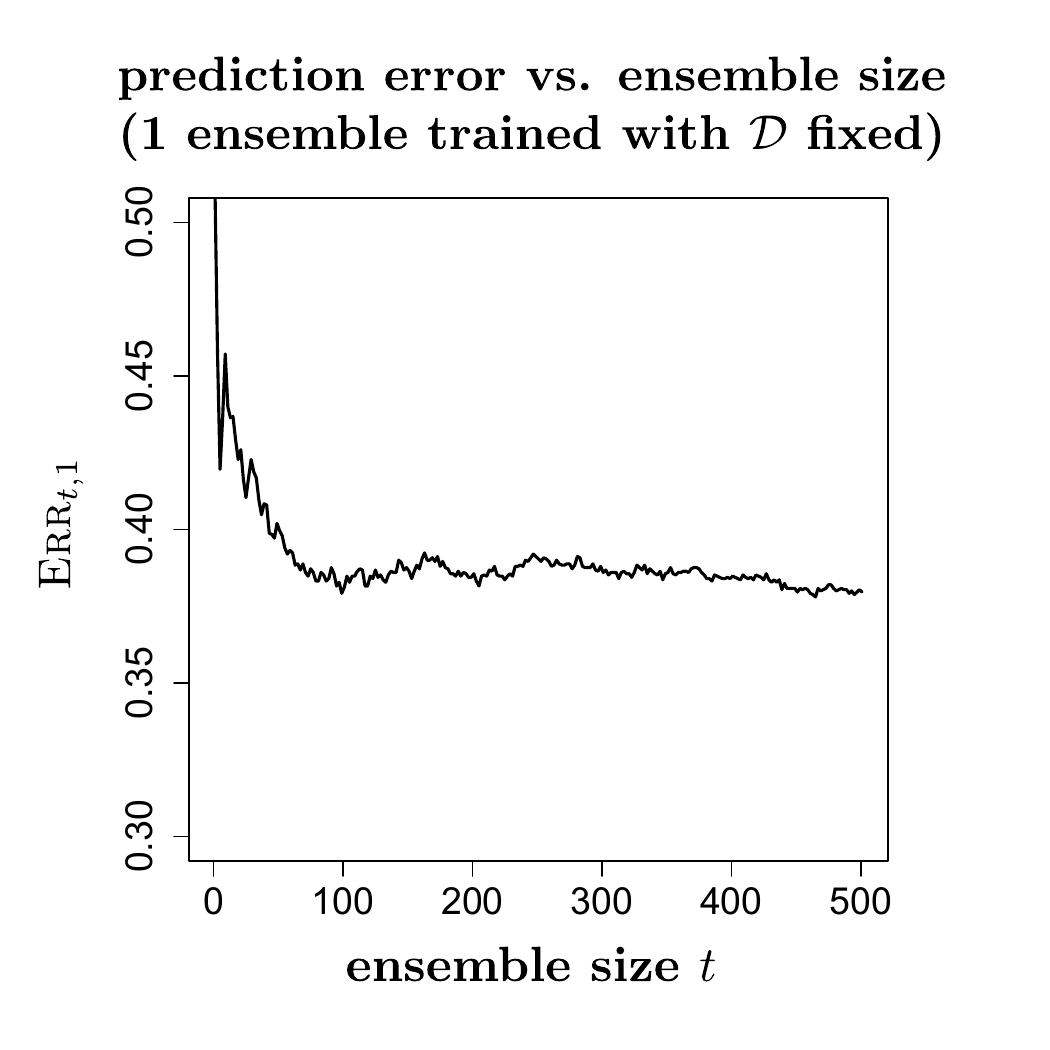}}
  \ \ \ \ 
  {\includegraphics[angle=0,
  width=0.45\linewidth]{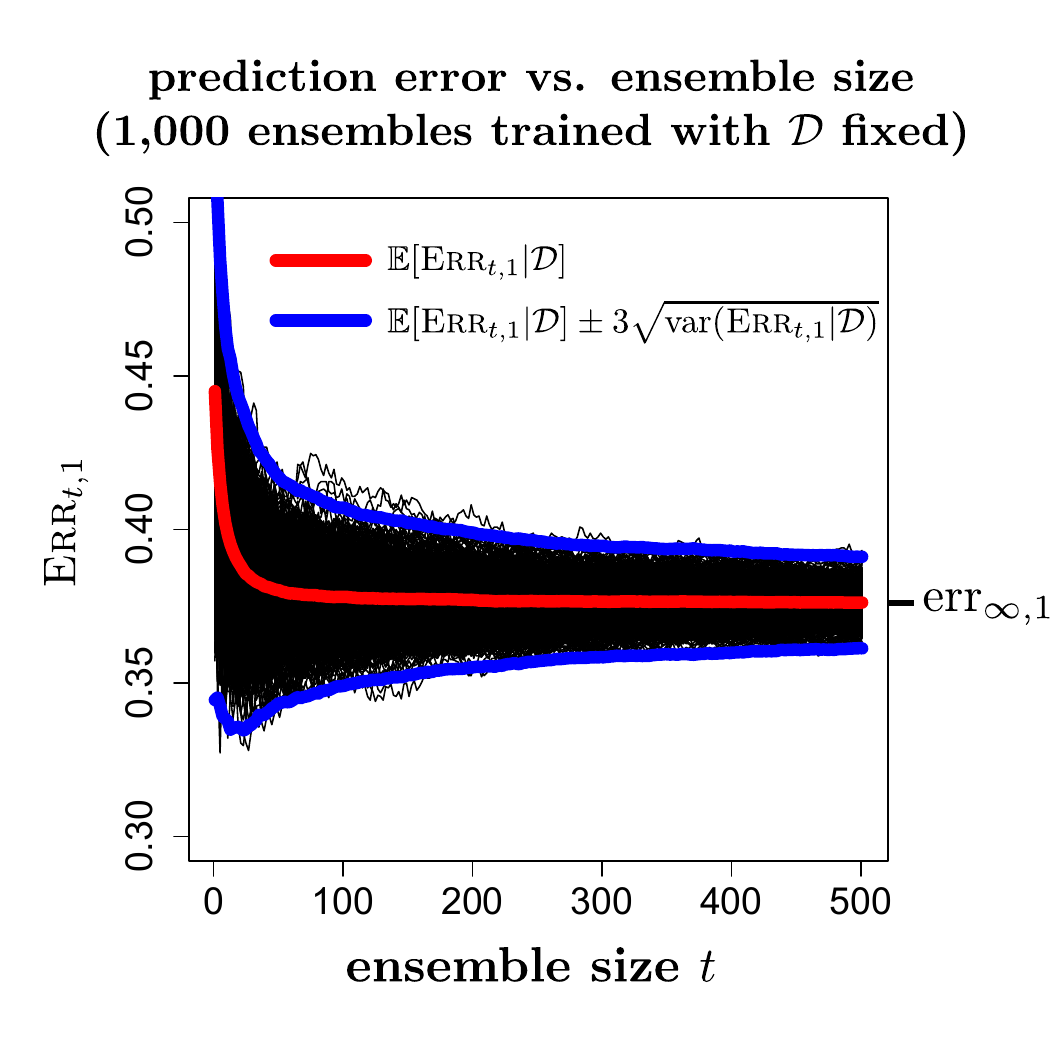}}
\caption{
 (Left panel): The fluctuations of $\Err_{t,1}$ for a single run of random forests on a fixed set $\D$. (Right panel): The fluctuations of $\Err_{t,1}$ for 1,000 runs of random forests on the same set $\D$.
  The value $\err_{\infty,1}$ represents the limit of $\Err_{t,1}$, as achieved by an infinite ensemble on $\D$.
 }
  \label{fig:intro}
 \end{figure}

To illustrate the randomness of $\Err_{t,l}$ when $\D$ is held fixed, the left panel of  Figure~\ref{fig:intro}  shows how $\Err_{t,l}$ fluctuates as an ensemble of 500 decision trees is trained by random forests.  
 (We refer to the curve in the left panel as a ``sample path''.) If this process is repeated many times,
 and random forests is used to train 1,000 ensembles (each with $t=500$) on the same dataset $\D$, then a pattern emerges:
We see 1,000 overlapping sample paths of $\Err_{t,l}$, in the right panel of Figure~\ref{fig:intro}.

Conceptually, the right panel gives some helpful insight into the algorithmic convergence of $\Err_{t,l}$. The main point to notice is that the convergence is well summarized by the mean and variance of $\Err_{t,l}$, conditionally on $\D$. Specifically, the red curve results from averaging the sample paths at each value of $t$, and hence represents $\E[\Err_{t,l}\big| \D]$, where the expectation is over the variables $\boldsymbol \xi_t=(\xi_1,\dots,\xi_t)$. Similarly, the variance of $\Err_{t,l}$ among the sample paths represents the \emph{algorithmic variance},
\begin{equation*}
\var(\Err_{t,l}\big| \D):=\E[\Err_{t,l}^2|\D]-\big(\E[\Err_{t,l}| \D]\big)^2,
\end{equation*}
which is the variance of $\Err_{t,l}$ due only to the training algorithm.
In this notation, the blue curves in Figure~\ref{fig:intro} are obtained by adding and subtracting $3\sqrt{\var(\Err_{t,l}|\D)}$ from the mean $\E[\Err_{t,l}|\D]$.

\paragraph{Problem formulation} \  As the ensemble size becomes large ($t\to\infty$),
 the prediction error $\Err_{t,l}$ typically converges in probability to a limiting value, denoted $\err_{\infty,l}$.
When judging the performance of bagging, random forests, or other methods satisfying \textbf{A1}, the value $\err_{\infty,l}$ plays a special role, since it is generally viewed as an ensemble's ideal class-wise error rate.
In this way, if $\Err_{t,l}$ is  sufficiently close to $\err_{\infty,l}$, then the ensemble has reached algorithmic convergence. Likewise, the problem of interest is to find a method for bounding the (unknown) gap $\Err_{t,l}-\err_{\infty,l}$ as a function of $t$.

Since $\Err_{t,l}$ is random, it is natural to measure the gap is in terms of its mean-squared value, which has the following bias-variance decomposition,
\begin{equation*}
\E\big[ (\Err_{t,l}-\err_{\infty,l})^2 \, \big| \, \D\big] =  \var(\Err_{t,l}|\D)+\big(\E[\Err_{t,l}|\D]-\err_{\infty,l}\big)^2. 
\end{equation*}
When comparing these bias and variance terms, 
it is important to note from the right panel of Fig.~\ref{fig:intro} that as long as $t$ is of modest size, then the bias $\E[\Err_{t,l}|\D]-\err_{\infty,l}$ is negligible compared to the standard deviation $\sqrt{\var(\Err_{t,l}|\D)}$. (Note that the bias is the difference between the red curve and $\err_{\infty,l}$, whereas the standard deviation controls the difference between the blue curves.)  Hence, the plot indicates that the size of the gap $\Err_{t,l}-\err_{\infty,l}$ is primarily governed by the standard deviation, and this is supported more generally by our theoretical results. Indeed, under certain assumptions, Lemma 1 shows that the bias is of order $\mathcal{O}(\ts\frac 1t)$, whereas Theorem 1 and Theorem 2 show that the standard deviation can be of order $\mathcal{O}(\ts\frac{1}{\sqrt{t}})$.

Based on these considerations, the problem of bounding the gap $\Err_{t,l}-\err_{\infty,l}$ can be reframed as the problem of bounding the variance $\var(\Err_{t,l}|\D)$. However, there is still a significant obstacle, because the variance describes how $\Err_{t,l}$ varies over \emph{repeated runs} of the ensemble method, while the user only has access to information from a \emph{single} ensemble. As a solution, our work shows that $\var(\Err_{t,l}|\D)$ is asymptotically bounded by a certain parameter that can be estimated effectively with a single ensemble.

\subsection{Contributions and related work}\label{sec:related}

\paragraph{Results on majority voting}
Our first main result is an upper bound on the algorithmic variance, which takes the form
$\var(\Err_{t,l}\big| \D) \leq \ts\frac{1}{4t} [f_l(\ts\frac{1}{2})]^2+o(\frac{1}{t}),$
where $f_l$ is a density function related to the ensemble method. The bound is presented in Theorem~\ref{VARTHM}, and the density $f_l$ will be defined there. As a complement to this result, we show in Theorem~\ref{ATTAINTHM} that the bound is sharp --- in the sense that it is attained by a specific family of randomized classifiers. In addition, we show in Corollary~\ref{EXCHCOR} that our analysis of $\Err_{t,l}$ can be extended to analyze stochastic processes beyond the context of classification. Specifically, our results can describe the running majority vote of general exchangeable Bernoulli sequences.

\paragraph{Methodology} \ \   To use the variance bound in practice, it is necessary to estimate the parameter $f_l(\ts\frac{1}{2})$ from a single run of the ensemble method. For this purpose, we propose two different estimators --- one based on a hold-out set, $\hat{f}_{l,\textsc{h}}(\ts\frac{1}{2})$, and another based on ``out-of-bag'' (OOB) samples, $\hat{f}_{l,\textsc{o}}(\ts\frac{1}{2})$. (See Section~\ref{sec:est} for a description of OOB samples.)  With regard to finite-sample performance, our experiments show that the resulting estimated bounds are tight enough to be meaningful diagnostics for convergence. 

Another important feature of the estimated bounds is that they are \emph{very inexpensive} to compute. In particular, once an ensemble has been tested with hold-out or OOB samples, each of the quantities $\hat{f}_{l,\textsc{h}}(\ts\frac{1}{2})$ and $\hat{f}_{l,\textsc{o}}(\ts\frac{1}{2})$ can be obtained with a \emph{single evaluation} of a kernel density estimator.
In this way, the bounds are an extra source of information that comes almost ``for free'' with the ensemble.

 \paragraph{MSE bound for the hold-out estimator} To analyze the performance of the hold-out estimator $\hat{f}_{l,\textsc{h}}(\ts\frac{1}{2})$, we derive an upper bound on its mean-squared error (MSE).
 An interesting aspect of the MSE bound is that it exhibits an ``elbow phenomenon,'' where the rate depends on the relative sizes of the hold-out set and $t$. 
Furthermore, when $t$ is sufficiently large, the bound turns out to match optimal non-parametric rates for estimating a density at a point.

\paragraph{Related work} \ \ For the most part, the literature on the problem of choosing the size of majority voting ensembles has focused on empirical approaches, e.g.~\cite{latinne,comet,swiss,oshiro}. Nevertheless, there have been a handful of previous works offering theoretical convergence analyses of $\Err_{t,l}$ and related quantities, as reviewed below.

The convergence of the expected error rate $\E[\Err_{t,l}\big| \D]$ seems to have been first studied in the paper~\cite{NgJordan}. Specifically, if the class proportions are denoted by $\pi_l=\P(Y=l)$, then the authors show that the total prediction error $\pi_0 \E[\Err_{t,0}\big| \D]+\pi_1 \E[\Err_{t,1}\big| \D]$ converges to its limit at the rate $1/t$, but the limiting constant is left unspecified. The first paper to determine the limiting constant was the 2013 preprint version of the current work~\cite{lopes2013}, showing that $\E[\Err_{t,l}\big| \D] \, = \, \text{{err}}_{\infty,l} +\ts\frac{1}{8t} f_l'(\ts\frac{1}{2})+o(\ts\frac{1}{t})$, as in line~\eqref{lemmathird} of Lemma~\ref{EXPECLEMMA} below. Later on, the 2015 preprint~\cite{cannings2015} obtained a corresponding formula under somewhat weaker regularity condition on $f_l$, and a relaxed form of the majority voting rule (namely $1\{\bar Q_t(X)\geq \alpha\}$, where $\alpha$ need not be $1/2$). Next, the 2016 preprint version of the current work~\cite{lopes2016} went beyond the expected value to establish an upper bound on $\var(\Err_{t,l}|\D)$, and demonstrated it to be sharp --- as in Theorems~\ref{VARTHM} and~\ref{ATTAINTHM} below. More recently, the variance bound was cited by the 2017 paper~\cite{cannings2017} (the journal version of~\cite{cannings2015}), and a similar result was included in its supplementary material. However, at a technical level, the variance bounds differ somewhat.  For instance, if $F_l$ denotes the cdf of $f_l$, then the paper~\cite{cannings2017} requires $F_l$ to be twice differentiable at $1/2$, whereas our current result allows $f_l$ to be non-differentiable at 1/2, provided that it is Lipschitz near $1/2$.  On the other hand, the result in~\cite{cannings2017} allows for the relaxed form of majority voting mentioned earlier.
Lastly, in terms of proofs, the result in~\cite{cannings2017} is derived from extensive calculations based on Edgeworth expansions, whereas the current proof is able to avoid Edgeworth expansions by using empirical process techniques. 

In another direction, the paper~\cite{hernandez2013} studied algorithmic convergence in terms of a different criterion, namely the ``disagreement probability'' $\delta_t:=\P(M_t(X)\neq M_{\infty}(X)\big| \D)$, where $M_{\infty}(X)$ is the infinite ensemble analogue of the majority vote $M_t(X)$. In that work, an informal derivation is given to show that $\delta_t$ is of order $\mathcal{O}(\ts\frac{1}{\sqrt{t}})$. However, as this relates to the current paper, the analysis of $\delta_t$ is a distinct problem, insofar as $\delta_t$ does not seem to provide a way to control the gap  $\Err_{t,l}-\err_{\infty,l}$,  or the variance $\var(\Err_{t,l}|\D)$. Also, estimation guarantees for $\delta_t$ have not previously been established. 
Nevertheless, it may be possible to obtain such guarantees using our results, since the leading term in the derivation for $\delta_t$ depends on the parameters $f_0(\ts\frac{1}{2})$ and $f_1(\ts\frac{1}{2})$.

Finally, with regard to methods that are supported by theoretical guarantees, an alternative approach to the current one is considered in the paper~\cite{lopes2019}. In that work, a bootstrap method is proposed to directly estimate $\var(\Err_{t,l}\big| \D)$. Two advantages of that approach are that it avoids the conservativeness of upper bounds, and it applies to the multi-class setting. However, the main advantage of the current method is that it has far lower computational cost. Indeed, if we suppose that an ensemble has already been tested with OOB samples (as is typically done by default), then the cost of the bootstrap method is at least of order $B\cdot t\cdot n$, where $B$ is the number of bootstrap samples. By contrast, the cost of evaluating the kernel density estimator is only $\mathcal{O}(n)$. Hence, if the user has a limited computational budget, then the current method allows the user to devote much more computation to other purposes --- such as training more classifiers, or optimizing an ensemble's tuning parameters.

\paragraph{Outline} \ \ 
In Section~\ref{sec:main}, we state some theoretical results on majority voting, which motivates the estimation method and guarantees given in Section~\ref{sec:est}. Later on, in Section~\ref{sec:expt}, we present some experiments to evaluate the estimation method.
In the supplementary material, we provide proofs for all theoretical results, as well as empirical validation of a technical assumption~(\textbf{A2}). 
\paragraph{Notation and terminology}\ \  If $U$ and $V$ are generic random variables, or sets of random variables, we write $\mathcal{L}(U|V)$ to refer to the conditional distribution of $U$ given $V$. For a function $g:[0,1]\to \R$, we say that $g$ is Lipschitz in a neighborhood of $1/2$ if there are positive constants $\kappa$ and $\delta$, such that $|g(s)-g(s')|\leq \kappa |s-s'|$ for all $s,s'\in[1/2-\delta,1/2+\delta]$.

\section{Theoretical results on majority voting }\label{sec:main}
In order to state our theoretical results, define the function $\theta: \mathcal{X}\to [0,1]$ according to
\begin{equation*}
\theta(x):=\E[Q_1(x)\big| \D],
\end{equation*}
where $x\in\mathcal{X}$ and $\D$ are held fixed, and the expectation is over the randomizing parameter $\xi_1$. When a random test point $X$ is plugged into $\theta$, we obtain a random variable $\theta(X)$ in the interval [0,1], and this variable will play an important role in our analysis. In particular, we will use the following technical assumption.

\begin{shortassumption}\label{density}
\emph{For each $l\in\{0,1\}$, the distribution $\mathcal{L}(\theta(X)\big| \D,Y=l)$ has a density $f_l$ on the interval $[0,1]$ that is Lipschitz in a neighborhood of $1/2$.}
\end{shortassumption}
 
 To provide some intuition for this assumption, first note that the majority vote of an infinite ensemble assigns a point $x\in\mathcal{X}$ to class 1 if and only if $\theta(x)\geq 1/2$.
 Consequently, the set of points $\mathcal{B}:=\{x\in\mathcal{X}: \theta(x)=1/2\}$ can be viewed as an ``asymptotic decision boundary'' in the space $\mathcal{X}$. As this relates to assumption \textbf{A2}, note that if the density $f_l$ has a ``spike'' at 1/2, then this means that a test point $X$ from class $l$ can fall exactly on the boundary $\mathcal{B}$. From an intuitive standpoint, this situation should make it difficult for the ensemble to reach ``consensus'' on test points, because some fraction of them will be  completely ambiguous. On the other hand, if $f_l$ is continuous in a neighborhood of 1/2, then a test point will have zero probability of falling exactly on $\mathcal{B}$, which should make it easier for the ensemble to reach consensus. In this way, the regularity of $f_l$ near 1/2 is related to the rate of convergence. 
 
With regard to the existence of the density $f_l$ in assumption~\textbf{A2}, it is possible to show that when the feature space is Euclidean, say $\mathcal{X}=\R^p$, and when the distribution $\mathcal{L}(X|Y=l)$ has a density on $\R^p$, then $f_l$ will exist as long as the function $\theta$ is smooth. (We refer to the books~\cite[Thm.~10.4,\, Thm.~10.6]{simon} and~\cite[p.345]{cucker} for details.)
Also, in the cases of bagging and random forests, it has been argued that ``bootstrap averaging'' has a smoothing effect on ``rough'' functions (such as decision trees), and consequently, the function $\theta$ can be smooth in an approximate sense~\cite{buhlmannyu,bujastuetzle2006,efron2014}. Moreover, in Appendix B, we provide empirical examples showing that for random forests, the distribution $\mathcal{L}(\theta(X)|\D,Y=l)$ is well approximated by distributions that satisfy \textbf{A2}. (Similar examples involving other datasets have also been given in~\cite{lopes2019}.)
\subsection{Variance bound and expectation formula}

Before stating the bound on algorithmic variance, we give a second-order formula for $\E[\Err_{t,l}\big| \D]$, where the random variable $\Err_{t,l}$ is as defined in~\eqref{errdef}. To make the formula easier to interpret, it is convenient to formally define the limiting error rates as
$\err_{\infty,0} := 1-F_0(\ts\frac{1}{2})$ and $\err_{\infty,1}:=F_1(\ts\frac{1}{2})$,
where $F_l$ denotes the c.d.f. associated with the density $f_l$.

\begin{lemma}\label{EXPECLEMMA}Let $l\in\{0,1\}$ and suppose \textbf{A1} and \textbf{A2} hold. Then,  as $t\to\infty$ along odd integers,
\begin{equation}\label{lemmasecond}
\E[\Err_{t,l}\big| \D] \, = \, \text{\emph{err}}_{\infty,l} \, + \, \mathcal{O}(\ts\frac{1}{t}).
\end{equation}
Furthermore, if $f_l$ is also differentiable at 1/2, then as $t\to\infty$ along odd integers,
\begin{equation}\label{lemmathird}
\E[\Err_{t,l}\big| \D] \, = \, \text{\emph{err}}_{\infty,l} +\ts\frac{1}{8t} f_l'(\ts\frac{1}{2})+o(\ts\frac{1}{t}).
\end{equation}
\end{lemma}

\paragraph{Remarks}  
In terms of algorithmic convergence, this formula is useful because it clarifies the relative importance of the bias $\big(\E[\Err_{t,l}|\D] -\err_{\infty,l}\big)$ and the standard deviation $\sqrt{\var(\Err_{t,l}|\D)}$. The fact that the standard deviation is at most of order $1/\sqrt{t}$ is the content of the following bound.

\begin{theorem}\label{VARTHM}
Let $l\in\{0,1\}$, and suppose~\textbf{A1} and~\textbf{A2} hold. Then, as $t\to\infty$ along odd integers, 
\begin{equation}\label{varbound}
\text{\emph{var}}(\Err_{t,l}\big| \D) \leq \ts\frac{1}{4t} [f_l(\ts\frac{1}{2})]^2+o(\frac{1}{t}).
\end{equation}
\end{theorem}

\paragraph{Remarks} 
To interpret the role of $f_l(\ts\frac{1}{2})$, recall from our earlier discussion that this parameter measures the density of ambiguous test points at the decision boundary. Hence, larger values of $f_l(\ts\frac{1}{2})$, should make it harder for the ensemble to reach consensus, which intuitively corresponds to higher algorithmic variance.

A priori, one might imagine that the algorithmic variance $\var(\Err_{t,l}|\D)$ could depend on many characteristics of the test point distribution and the ensemble method. From this perspective, the bound~\eqref{varbound} has a surprisingly simple form. Even so, our next result shows that the bound cannot be improved in general (with respect to ensemble methods satisfying \textbf{A1} and \textbf{A2}).

\subsection{Attaining the variance bound}\label{sec:attain}

We now aim to show that the variance bound in Theorem~\ref{VARTHM} is attained by a specific family of classifiers. A notable feature of this construction is that it does not require ``pathological'' choices for the ensemble method. In fact, starting from \emph{any} choice of $\{Q_i\}$ that satisfies the conditions of Theorem~\ref{VARTHM}, it is possible to construct a related ensemble $\{Q_i^{\circ}\}$ that attains the bound. Furthermore, the ensembles $\{Q_i\}$ and $\{Q_i^{\circ}\}$ will turn out to have the same  class-wise error rates on average.

To proceed with the construction, let $\theta(x)=\E[Q_1(x)|\D]$ as before, and let $U_1,U_2,\dots$ be an i.i.d.~sequence of Uniform[0,1] variables 
(independent of $\D$ and $(X,Y)$). 
Next, for each $i=1,2,\dots,$  define the random classifier function $Q^{\circ}_i:\mathcal{X}\to \{0,1\}$, according to
\begin{equation*}
Q_i^{\circ}(x):=1\{\theta(x)\geq U_i\}.
\end{equation*}
As a way of making sense of this definition, recall that when $t=\infty$, the majority vote of the  ensemble $\{Q_i\}$ is given by the indicator $1\{\theta(x)\geq 1/2\}$.  Hence, the classifier $Q_i^{\circ}$ can be viewed a randomized version of the asymptotic majority vote, since the variable $U_i$ plays the role of $\xi_i$, and can be interpreted as a ``random threshold'' whose expected value is $1/2$.

By analogy with the original ensemble, we define the class-wise error rates of the new ensemble as 
$\Err^{\circ}_{t,l}:= \ts\int_{\mathcal{X}} 1\{|\bar{Q}^{\circ}_t(x)-l|\geq  \ts \frac{1}{2}\} d\mu_l(x),$
where again $\mu_l=\mathcal{L}(X|Y=l)$.
The following theorem shows that the ensemble $\{Q_i^{\circ}\}$ attains the highest possible algorithmic variance as $t\to\infty$ .

\begin{theorem}\label{ATTAINTHM} Let $l\in\{0,1\}$, and suppose  \textbf{A1} and \textbf{A2} hold. Then as $t\to\infty$ along odd integers,
\begin{equation*}
 \text{\emph{var}}(\Err^{\circ}_{t,l}\big| \D) = \ts\frac{1}{4t} [f_l(\ts\frac{1}{2})]^2+o(\frac{1}{t}).
\end{equation*}
\end{theorem}
\paragraph{Remarks} \  Regarding the relative performance of the two ensembles  $\{Q_i\}$ and $\{Q_i^{\circ}\}$, it is interesting to note that their class-wise error rates are equal on average,~$\E[\Err_{t,l}\big| \D]=\E[\Err_{t,l}^{\circ}\big| \D]$. Indeed, this can be checked by using the fact that for any fixed $x\in\mathcal{X}$, the binary sequences $\{Q_i(x)\}$ and $\{Q_i^{\circ}(x)\}$ are both i.i.d.~Bernoulli$(\theta(x))$, conditionally on $\D$.

 Using a version of Slepian's lemma, a fairly simple heuristic argument can be given to suggest why the ensemble $\{Q_i^{\circ}\}$ attains the variance bound.  The relevant version of this fact states that if $(Z_1,Z_2)$ is a centered bivariate normal random vector, with $\text{var}(Z_1)=\text{var}(Z_2)=1$ and $\text{cor}(Z_1,Z_2)=\rho$, then for any numbers $a,b\in\R$, the ``corner probability'' $\P(Z_1\leq a, Z_2\leq b)$
 is a non-decreasing function of $\rho$~\cite{slepian1962,sidak1968}. To apply this fact, first recall that the error rates $\Err_{t,l}$ and $\Err_{t,l}^{\circ}$ are equal on average, and so we may compare their variances by comparing their second moments. Some elementary manipulation of the definition~\eqref{errdef} gives the expression
\begin{equation*}
\E[\Err_{t,1}^2\big| \D]=\int_{\mathcal{X}}\int_{\mathcal{X}} \P\Big((Z_t(x),Z_t(x'))\in \mathcal{C}_t(x,x')\Big| \D \Big)d\mu_1(x)d\mu_1(x'),
\end{equation*}
where we define the random variable
\begin{equation*}
Z_t(x):=\sqrt{t}(\bar{Q}_t(x)-\theta(x)),
\end{equation*}
and the two-dimensional corner set 
\begin{equation*}
\mathcal{C}_t(x,x'):=(-\infty,c_t(x)]\times (-\infty,c_t(x')] \text{ \ \ \ with \ \ \ } c_t(x):=\sqrt{t}(1/2-\theta(x)).
\end{equation*}
As $t\to\infty$, the vector $(Z_t(x),Z_t(x'))$ approaches a centered bivariate Gaussian distribution.
Therefore, Slepian's lemma suggests that the second moment $\E[\Err_{t,1}^2\big| \D]$ can be bounded asymptotically by replacing $\{Q_i\}$ with $\{Q_i^{\circ}\}$, and then checking that for each pair $(x,x')$, the variables $Q_i^{\circ}(x)$ and $Q_i^{\circ}(x')$ are maximally correlated. The latter step works out easily because 
\begin{equation*}
\begin{split}
\E[Q_i(x)Q_i(x')\big| \D] &\leq \min\big\{\E[Q_i(x)\big| \D]\, ,\,\E[Q_i(x')\big| \D]\big\},\\[0.2cm]
&=\min\{\theta(x),\theta(x')\},\\[0.2cm]
&=\E[Q_i^{\circ}(x)Q_i^{\circ}(x')\big| \D].
\end{split}
\end{equation*}
Although this informal reasoning leads to the correct conclusion, the formal proof in the supplement will approach the problem differently in order to avoid some technical complications.

\subsection{A corollary for exchangeable Bernoulli sequences}

Exchangeable stochastic processes are a fundamental topic in probability and statistics, and in this subsection, we take a short sidebar to explain how our formula for $\E[\Err_{t,l}\big| \D]$ can be expressed in the language of exchangeability. The basic link between exchangeability and ensemble methods occurs through de Finetti's theorem~\cite[Ch.~1.4]{schervish2012}, which we now briefly review.

An infinite sequence of random variables $V_1,V_2,\dots$ is said to be \emph{exchangeable} if the joint distribution of any finite sub-collection is invariant under permutation. That is, $\mathcal{L}(V_{i_1},\dots,V_{i_k}) = \mathcal{L}(V_{\tau(i_1)},\dots,V_{\tau(i_k)})$, for all positive integers $k$, and all permutations $\tau$ on $k$ letters. In the special case that each $V_i$ is a Bernoulli random variable, de Finetti's theorem states that the sequence $V_1,V_2,\dots$ is exchangeable if and only if there is a random variable $\Theta$ in the unit interval $[0,1]$, such that, conditionally on $\Theta=\vartheta$, the variables $V_1,V_2,\dots$ are i.i.d.~Bernoulli$(\vartheta)$~\cite[Thm.~1.47]{schervish2012}. As a matter of terminology, the distribution of the random variable $\Theta$ is called the \emph{mixture distribution} associated with the sequence $\{V_i\}$.

If we define the running majority vote of an exchangeable Bernoulli sequence as the indicator
  \begin{equation*}
   {\texttt{\large{M}}}_t:=1\{\ts\frac{1}{t}\sum_{i=1}^t V_i \geq \frac{1}{2}\},
 \end{equation*}
then the following corollary shows that the expectation $\E[   {\texttt{\large{M}}}_t]$ obeys a second order formula analogous to the one in Lemma~\ref{EXPECLEMMA}.
\begin{corollary}\label{EXCHCOR}
Let $V_1,V_2,\dots$ be an infinite exchangeable Bernoulli sequence whose mixture distribution function is denoted by $F$.  Suppose the function $F$ has a density $f$ on $[0,1]$ that is Lipschitz in a neighborhood of $1/2$. Then as $t\to\infty$ along odd integers,
\begin{equation}\label{main}
1-\E[   {\texttt{\large{\emph{M}}}}_t] = F(\ts\frac{1}{2})+\mathcal{O}(\ts\frac{1}{t}).
\end{equation}
Furthermore, if $f$ is also differentiable at $1/2$, then
\begin{equation}\label{main}
1-\E[   {\texttt{\large{\emph{M}}}}_t] = F(\ts\frac{1}{2})+\ts\frac{1}{8t}f'(\ts\frac{1}{2})+o(\ts\frac{1}{t}).
\end{equation}
\end{corollary}
\paragraph{Remarks}  To see the connection with an ensemble of classifiers, de Finetti's theorem implies that each random variable $V_i$ can be considered in terms of a random binary function \mbox{$W_i:[0,1]\to \{0,1\}$} such that $V_i=W_i(\Theta)$, and for a fixed value $\Theta=\vartheta$, the random variables $W_1(\vartheta),W_2(\vartheta),\dots$ are i.i.d.~Bernoulli$(\vartheta)$. In other words, the functions $\{W_i\}$ play the role of the classifiers $\{Q_i\}$, and the mixing variable $\Theta$ plays the role of $\theta(X)$. Once this translation has been made, the proof of Lemma~\ref{EXPECLEMMA} carries over directly to Corollary~\ref{EXCHCOR}.

\section{Methodology and guarantees}\label{sec:est}
In this section, we present two methods for estimating the parameter $f_l(\ts\frac{1}{2})$, as well as a consistency result in Theorem~\ref{MSETHM}.

\subsection{Estimation with a hold-out set}\label{sec:heldout}
If it were possible to obtain samples directly from the density $f_l$, say $\Theta_1,\dots,\Theta_{m_l}$, then a natural approach to to estimating $f_l(\ts\frac{1}{2})$ would proceed via a kernel density estimator of the form
\begin{equation*}
 \ts\frac{1}{m_lh}\displaystyle\sum_{j=1}^{m_l} K\Big(\ts\frac{1/2-\Theta_j}{h}\Big),
\end{equation*}
where $K:\R\to \R$ is a kernel function satisfying $\int_{\R}K(s)ds=1$, and the number $h>0$ is a bandwidth parameter. However, the main difficulty we face is that direct samples from $f_l$ are unavailable in practice. Instead, we propose to construct ``noisy samples'' from $f_l$ along the following lines.

To proceed, suppose that a set of i.i.d.~samples $\tilde{X}_{1,l},\dots,\tilde{X}_{m_l,l}\sim \mathcal{L}(X|Y=l)$ from class $l$ have been held out. (In particular, these hold-out samples  are assumed to be independent of the training set $\D$, the ensemble $\{Q_i\}$, and the test point $(X,Y)$.) If the function $\theta$ were known exactly, the hold-out samples could be plugged into $\theta$ to create i.i.d.~samples $\theta(\tilde{X}_{1,l}),\dots,\theta(\tilde{X}_{m_l,l})$ from the distribution $f_l$. So, using the fact the averaged classifier $\bar{Q}_t$ approximates the function $\theta$ as $t\to\infty$, we may regard the observable values $\bar{Q}_t(\tilde{X}_{1,l}),\dots,\bar{Q}_t(\tilde{X}_{m_l,l})$ as noisy samples from $f_l$. Next, if we let the random variable $\ve_{j,l}$, for $j\in\{1,\dots,m_l\}$ be defined by
\begin{equation}\label{noisemodel}
\bar{Q}_t(\tilde{X}_{j,l}) = \theta(\tilde{X}_{j,l})+\ve_{j,l},
\end{equation}
then $\ve_{j,l}$ can be interpreted as noise with mean zero. From a deconvolution perspective, this model is challenging, since $\var(\ve_{j,l})$ is \emph{unknown}, and  also, $\ve_{j,l}$ is \emph{not independent} of $\theta(X_{j,l})$. Nevertheless, it is plausible that the estimation of $f_l(\ts\frac{1}{2})$ is still tractable, since the following bound shows that the noise becomes small as $t$ increases,
\begin{equation*}
\begin{split}
\var(\ve_{j,l}\big|\D)
&=\E[\ts\frac{1}{t}\theta(\tilde{X}_{j,l})(1-\theta(\tilde{X}_{j,l}))|\D]\\[0.2cm]
&\leq \ts\frac{1}{4t},
\end{split}
\end{equation*}
where the first line follows from the law of total variance.
Consequently, we propose to estimate $f_l(\ts\frac{1}{2})$ by directly applying a kernel density estimator to the values $\bar{Q}_t(\tilde{X}_{j,l})$, with the hold-out estimator defined as

\begin{equation*}
\hat{f}_{l,\textsc{h}}(\ts\frac{1}{2}) := \ts\frac{1}{m_l h}\displaystyle\sum_{j=1}^{m_l} K\Big(\ts\frac{1/2-\bar{Q}_t(\tilde{X}_{j,l})}{h}\Big),
\end{equation*}
For theoretical convenience, we will only analyze the ``rectangular kernel''
\mbox{$K(s):=\ts\frac{1}{2}1\{-1\leq s\leq 1\},$}
but our method can be applied to any choice of kernel in practice. In the next subsection, we will specify the size of the bandwidth as an explicit function of $m_l$ and $t$.

In terms of computation, the bulk of the cost to calculate $\hat{f}_{l,\textsc{h}}(1/2)$ comes from obtaining the values $\bar{Q}_t(\tilde{X}_{j,l})$. Often, these values are computed anyway when estimating the error rate $\Err_{t,l}$ from a hold-out set. Hence, the extra information provided by the estimator $\hat{f}_{l,\textsc{h}}(1/2)$ comes at a very small added cost. Furthermore, the same computational benefit holds for the ``OOB estimator'' proposed in Section~\ref{sec:oob}.

\subsection{An MSE bound for the hold-out estimator}
 \ To measure the accuracy of the hold-out estimator, we consider the mean-squared error,
$$\textsc{mse}(\hat{f}_{l,\textsc{h}}(\ts\frac{1}{2})) := \E\big[(\hat{f}_{l,\textsc{h}}(\ts\frac{1}{2})-f_l(\ts\frac{1}{2}))^2\big| \D\big].
$$
Here, the expectation is over both the hold-out set $\tilde{X}_{1,l},\dots,\tilde{X}_{m_l,l}$, and the randomizing variables $\{\xi_i\}$. Although the conditioning on $\D$ may appear unusual in this definition of MSE, it is necessary because the parameter $f_l(1/2)$ is specific to the dataset $\D$ (since the function $\theta$ is). The following result gives a non-asymptotic bound on the MSE, which holds for fixed values of  $t$ and  $m_l$. 

\begin{theorem}\label{MSETHM}
Under  \textbf{A1} and \textbf{A2}, let $[1/2-\delta_l,1/2+\delta_l]$ denote a neighborhood on which $f_l$ is Lipschitz, with $\delta_l\in(0,1/2)$, and $l\in\{0,1\}$.
 Also, suppose  $\hat{f}_{l,\textsc{h}}(\ts\frac{1}{2})$ is computed with the rectangular kernel. Under these conditions, there are numbers $c_1,c_2>0$ not depending on $t$ or $m_l$,  such that if the bandwidth is set to $h=c_1(\min\{m_l,\sqrt{t}\})^{-1/3}$,   then the following bound holds as soon as  $h<\delta_l$,
\begin{equation}\label{mse1}
 \textsc{mse}(\hat{f}_{l,\textsc{h}}(\ts\frac{1}{2})) \ \leq \
c_2\big(\min\{m_l,\sqrt{t}\}\big)^{-2/3}.
\end{equation}
Furthermore, if the derivative $f_l'$ is also Lipschitz on $[1/2-\delta_l,1/2+\delta_l]$, and if the bandwidth is set to $h=c_1(\min\{m_l,\sqrt{t}\})^{-1/5}$, then the rate $-2/3$ above may be replaced with $-4/5$.
\end{theorem}

\paragraph{Remarks}  A notable aspect of the result is that the MSE bound has an intertwined dependence on the sample size $m_l$ and computational cost $t$. This connection is also interesting because it presents an ``elbow phenomenon'', where the bound qualitatively changes, depending on whether $m_l< \sqrt{t}$ or $m_l>\sqrt{t}$.

 With regard to minimax optimality, it is clear that the bound's dependence on $m_l$ cannot be improved --- provided that estimation is based on the values $\bar{Q}_t(\tilde{X}_{1,l}),\dots,\bar{Q}_t(\tilde{X}_{m_l,l})$. To see this, note that when $t=\infty$, the problem reduces to estimating $f_l(1/2)$ with \emph{noiseless} i.i.d.~samples $\theta(\tilde{X}_{1,l}),\dots,\theta(\tilde{X}_{m_l,l})\sim f_l$. In this case, it is well known that if $f_l$ is restricted to lie in a class of densities $g$ for which the $(\beta-1)$th derivative $g^{(\beta -1)}$ is Lipschitz, then the rate $m_l^{-2\beta/(2\beta+1)}$ is optimal~\cite{tsybakov}. On this point, it is somewhat surprising that the ``noiseless rate'' $m_l^{-2\beta/(2\beta+1)}$ ``kicks in'' as soon as $m_l<\sqrt{t}$, because for finite values of $t$, the estimator $\hat{f}_{l,\textsc{h}}(1/2)$ is built from noisy samples --- and in deconvolution problems, the optimal rates are typically slower~\cite[Theorem 2.9]{meister}.

 To some extent, this effect of attaining the noiseless rate for sufficiently large $t$ may be explained by the fact that the noise variance scales like $1/t$ in the model~\eqref{noisemodel}. However,  the overall situation is complicated by the fact that the variables $\ve_{j,t}$ and $\theta(\tilde{X}_{j,l})$ are not independent, and the noise variance unknown. In the deconvolution literature, a few other works have reported on a similar phenomenon of attaining ``fast'' convergence rates when the noise level is ``small'' in various senses~\cite{delaigle,hesse,JMVAdecon,meister2010}.
Nevertheless, the models in these works are not directly comparable with the model~\eqref{noisemodel}.  Likewise, we leave a more detailed analysis of the model~\eqref{noisemodel} for future work, since our main focus is on measuring the algorithmic convergence of randomized ensembles.

\subsection{Estimation with out-of-bag samples}\label{sec:oob}
To avoid the need for a hold-out set, the previous estimator can be modified to take advantage of ``out-of-bag'' (OOB) samples, which are a special feature of bagging and random forests. For a quick description of OOB samples, note that when bagging and random forests are implemented, each classifier $Q_i$ is trained on randomly selected set $\mathcal{D}_i^*$,  obtained from $\mathcal{D}$ by sampling with replacement. Due to this sampling mechanism, approximately $(1-1/n)^n\approx 37\%$ of the $n$ training samples in $\mathcal{D}$ are likely to be excluded from each $\mathcal{D}_i^*$ --- and these excluded (OOB) samples are useful because they serve as ``effective test samples'' for each classifier.

As a matter of terminology, if a training point $X_j$ is not included in $\D_i^*$, we will say that $X_j$ is out-of-bag for the classifier $Q_i$. 
Likewise, for each index $j\in\{1,\dots,n\}$, we define the set \mbox{$\textsc{oob}(j)\subset\{1,\dots,t\}$} to index the classifiers for which $X_j$ is out-of-bag. Hence, by fixing attention on a test point $X_j$, and then averaging over the values $Q_i(X_j)$ with $i\in\textsc{oob}(j)$, we can obtain an approximate sample from the distribution $f_l$. 

To finish carrying out this idea, let $j\in\{1,\dots,n\}$, and let the random function $\bar{Q}_{t,j}:\mathcal{X}\to [0,1]$ be defined by
\begin{equation*}
\bar{Q}_{t,j}(x):=\ts\frac{1}{|\textsc{oob}(j)|}\tsum_{i\in \textsc{oob}(j)}Q_i(x),
\end{equation*}
 where $|\textsc{oob}(j)|$ denotes the cardinality of $\textsc{oob}(j)$, and we put $\bar{Q}_{t,j}(x)=0$ in the rare case that $\textsc{oob}(j)$ is empty.\footnote{Note that $\textsc{oob}(j)$ is empty with probability $[1-(1-1/n)^{n}]^t \approx(0.63)^t$.} Next, if we let $\mathcal{J}_l\subset\{1,\dots,n\}$ index the training points from class $l$, then we  define the the OOB estimator for $f_l(\ts\frac{1}{2})$ as
 \begin{equation*}
 \hat{f}_{l,\textsc{o}}(\ts\frac{1}{2}) := \ts\frac{1}{|\mathcal{J}_l| h}\tsum_{j\in\mathcal{J}_l} K\Big(\ts\frac{\bar{Q}_{t,j}(X_{j})-1/2}{h}\Big),
\end{equation*}
for a given choice of kernel $K$ and bandwidth $h$. Lastly, it is worth emphasizing that the values $\bar{Q}_{t,j}(X_{j})$ are often computed by default during a run of bagging or random forests, because these values are used to compute the ``out-of-bag error rate'', which is a standard alternative to a hold-out error estimate.

\section{Numerical experiments}\label{sec:expt}

\paragraph{Design of experiments}\label{sec:design}  The goal of the experiments is to look at how close the estimated bounds $\hat{f}_{l,\textsc{h}}(1/2)/(2\sqrt{t})$ and $\hat{f}_{l,\textsc{o}}(1/2)/(2\sqrt{t})$ are to the unknown quantity $\sqrt{\var(\Err_{t,l}|\mathcal{D})}$.
The experiments were based on classification problems associated with the following datasets from the UCI repository~\cite{uci}:  `abalone', `optical recognition of handwritten digits' (abbrev.~`digits'), `HIV-1 protease cleavage' (abbrev.~`HIV'),  `landsat satellite', `occupancy detection', and `spambase'. Additional details regarding data preparation are discussed at the end of the supplement.

Each full dataset of labeled examples was evenly split into a training set $\D$, and a separate ``ground truth'' set $\mathcal{D}_{\text{ground}}$, with nearly matching class proportions in $\D$ and $\mathcal{D}_{\text{ground}}$. The set $\mathcal{D}_{\text{ground}}$ was used for approximating the ground truth values of $\E[\Err_{t,l}|\mathcal{D}]$ and $\var(\Err_{t,l}|\D)$, and that is why a substantial amount of data was reserved for $\mathcal{D}_{\text{ground}}$. Also, a smaller set $\mathcal{D}_{\text{hold}}\subset \mathcal{D}_{\text{ground}}$ of size $|\mathcal{D}_{\text{hold}}| \approx 0.2 |\mathcal{D}|$ was used as the ``hold-out'' set  for computing the hold-out estimator $\hat{f}_{\textsc{h},l}(1/2)$. (The smaller size of $\mathcal{D}_{\text{hold}}$ was chosen to illustrate the performance of $\hat{f}_{\textsc{h},l}(1/2)$ from a limited hold-out set.)

 Using standard settings, the R package `randomForest'~\cite{randomForestsCitation} was used to train 5,000 ensembles on $\mathcal{D}$, with each containing $t=500$ classifiers (the default size). In turn, each ensemble was tested on $\mathcal{D}_{\text{ground}}$, producing 5,000 estimates of $\Err_{t,l}$ for each class $l\in\{0,1\}$. In the table below, the sample mean and standard deviation of these 5,000 estimates are reported as $\E[\Err_{t,l}|\mathcal{D}]$ and $\sqrt{\var(\Err_{t,l}|\mathcal{D})}$ respectively. These two values are viewed as ground truth, but of course, they are imperfect, since their quality is limited by the size of $\mathcal{D}_{\text{ground}}$ and monte-carlo error.

To estimate the theoretical bound $f_l(1/2)/(2\sqrt{t})$ on $\sqrt{\var(\Err_{t,l}|\D)}$, we implemented the hold-out and OOB methods from Section~\ref{sec:est}. Specifically, the estimators $\hat{f}_{l,\textsc{h}}(1/2)/(2\sqrt{t})$ and $\hat{f}_{l,\textsc{o}}(1/2)/(2\sqrt{t})$ were computed for each ensemble, giving 5,000 realizations of each estimator. Also, for each estimator, we used the rectangular kernel, and the R default bandwidth selection rule `nrd0'.
 In the table below, the sample average of the 5,000 realizations of $\hat{f}_{l,\textsc{h}}(1/2)/(2\sqrt{t})$ is referred to as the `hold-out bound', and the sample standard deviation is listed in parentheses. The results for the OOB estimator are reported similarly, under the name `OOB bound'.

\begin{table}[h!]
\centering
\scriptsize
\setlength\extrarowheight{3pt}
\begin{tabular}{lllllll}
\multicolumn{1}{l}{} &   \multicolumn{2}{c}{`abalone' } & \multicolumn{2}{c}{`digits'} & \multicolumn{2}{c}{`HIV'}\\ 
 \cline{2-7}
   & class 0 & class 1  & class 0 & class 1  & class 0 & class 1   \\
   \cline{1-7}
\multicolumn{1}{l}{$\E[\Err_{t,l}|\mathcal{D}]$ \ (\%)}  
   & 28.99  & 7.72 & 1.49  & 2.44 &  5.24 & 12.67    \\\cline{1-7}
\multicolumn{1}{l}{$\sqrt{\var(\Err_{t,l}|\mathcal{D})}$ \ (\%)   } 
& .37 & .20 & .15 & .16 & .15 & .48 \\\cline{1-7}
\multicolumn{1}{l}{hold-out bound  (\%)}    
& 1.94 (.05)  & .69 (.08) & .31 (.08) &  .68 (.09) & .95 (.08) & 1.53 (.10)                 \\ \cline{1-7}
\multicolumn{1}{l}{OOB bound  (\%)}
      &  1.50 (.05)  & .86 (.06) &  .34 (.06)  & .41 (.07) & .90 (.06)  & 1.75 (.09) \\ \cline{1-7}
\multicolumn{1}{l}{\# hold-out samples}
    & 145     & 273    & 280     & 283     & 527  & 133 \\ \cline{1-7}
\multicolumn{1}{l}{\# training samples}  
  & 683     & 1,406   & 1,423     & 1,387    & 2,598  & 697   \\ \wcline{1-7}
\end{tabular}
\setlength\extrarowheight{3pt}
\begin{tabular}{lllllll}
\multicolumn{1}{l}{} &   \multicolumn{2}{c}{`landsat satellite'} & \multicolumn{2}{c}{`occupancy detection' } & \multicolumn{2}{c}{`spambase'}\\
 \cline{2-7}
   & class 0 & class 1  & class 0 & class 1  & class 0 & class 1   \\
   \cline{1-7}
\multicolumn{1}{l}{$\E[\Err_{t,l}|\mathcal{D}]$ \ (\%)} 
& 2.50   & 7.05   & 5.41  & 29.87   & 2.89 & 10.02    \\\cline{1-7}
\multicolumn{1}{l}{$\sqrt{\var(\Err_{t,l}|\mathcal{D})}$ \ (\%)   }
& .11 & .13  & .06 & .11 & .10 & .23   \\\cline{1-7}
\multicolumn{1}{l}{hold-out bound  (\%)}
& .55 (.10)   & .17 (.11)  & .39 (.06) &  .47 (.02) & .24 (.09) & .52 (.07)
               \\ \cline{1-7}
\multicolumn{1}{l}{OOB bound  (\%)}
 &  .25 (.05)    & .51 (.10)   &  .42 (.04) & .48 (.02) & .36 (.07) &  .45 (.07)
         \\ \cline{1-7}
\multicolumn{1}{l}{\# hold-out samples}
 & 356  & 289      &  1,477   & 580  & 279  & 182        \\ \cline{1-7}
\multicolumn{1}{l}{\# training samples}
&  1,818   &  1,400    & 7,271   & 3,009     &  1,397  &  904   \\ \cline{1-7}
\end{tabular}\\
\caption{(Numerical results for hold-out and OOB bounds). The quantities $\Err_{t,l}|\mathcal{D})$ $\sqrt{\var(\Err_{t,l}|\mathcal{D})}$ refer to the mean and standard deviation of $\Err_{t,l}$, conditionally on a given dataset, over 5,000 repeated runs of the random forests algorithm, as described in the main text. Next, the quantities `hold-out bound' and `OOB bound' respectively refer to $\hat{f}_{l,\textsc{h}}(1/2)/(2\sqrt{t})$ and $\hat{f}_{l,\textsc{o}}(1/2)/(2\sqrt{t})$, and the corresponding entries in the table refer to the mean and standard deviations (in parentheses) over the 5,000 repeated runs. Lastly, all rows labeled with $(\%)$ are reported in percentage values, so that 1 corresponds to .01.}
\end{table}
\normalsize

\paragraph{Comments on results}

Before discussing the results in detail, it is worth clarifying  that even when the random forests method has poor accuracy as a classifier, it is still possible for the estimated bounds to serve their purpose. An example of this occurs in `occupancy detection', class 1, where the estimated bounds are reasonably tight, even though the expected error rate is almost 30\%.

There are several main conclusions we can draw from the table. The first point to note is that in all cases, the estimated  bounds are indeed larger than $\sqrt{\var(\Err_{t,l}|\D)}$. Although this is what we expect, it is not an entirely trivial property, because this could be violated if the estimators $\hat{f}_{l,\textsc{h}}(\ts\frac{1}{2})$  and $\hat{f}_{l,\textsc{o}}(\ts\frac{1}{2})$ are not sufficiently close to $f_l(\ts\frac{1}{2})$.
Second, the hold-out and OOB methods have mostly similar performance across the datasets. Consequently, the OOB method is likely to be preferred in practice, since it does not require a hold-out set.~Third, even though the bounds are  fairly conservative, they can still be tight enough to provide meaningful information about algorithmic convergence. For instance, in 10 of the 12 cases, both bounds are able to confirm that $\sqrt{\var(\Err_{t,l}|\D)}$ is less than 1\%. Also, the bounds are usually much smaller than $\E[\Err_{t,l}|\D]$, which is a notable point of reference, because it is natural to judge the fluctuations of $\Err_{t,l}$ in relation to the size of the error itself. (For instance, when the error rate is high, it may be reasonable to tolerate larger fluctuations.)

\section*{Appendices}  
  
 Appendix A includes all proofs. Appendix B discusses empirical validation of assumption \textbf{A2}. Appendix C discusses details regarding data preparation. 
\appendix

\section{Proofs}\label{app:proofs}

For simplicity, in most of the proofs, we will only treat the case of $l=1$, since the proofs for $l=0$ are essentially identical. To simplify notation, we will allow $\kappa>0$ to denote a constant not depending on $t$ or $m_l$, whose value may differ from line to line. Another piece of notation is that $\xrightarrow{ \ d \ }$ refers to convergence in distribution.  Recall also that $t$ is always odd. 

Results are proven in the same order as they appear in the main text. Specifically,  Lemma 1 is proven in~\ref{EXPECLEMMAproof}, Theorem 1 is proven in~\ref{VARTHMproof}, Theorem 2 is proven in~\ref{ATTAINTHMproof}, and Theorem 3 is proven in~\ref{MSETHMproof}

\subsection{Proof of Lemma~\ref{EXPECLEMMA} (expectation formula)}\label{EXPECLEMMAproof}
We first prove Lemma 1 in the case when $f_1$ is differentiable at 1/2, and then later in~\ref{nondiffcase}, we explain the small modification needed to handle the case when the differentiability condition does not hold.

To begin with some notation, define the random binary function $H_t:\mathcal{X}\to \{0,1\}$ according to
\begin{equation}\label{Htdef}
H_t(x):=1\{\bar{Q}_t(x)\leq \ts\frac{1}{2}\},
\end{equation}
which allows us to write
\begin{equation}\label{Hrep}
\Err_{t,1}=\ts\int_{\mathcal{X}} H_t(x)d\mu_1(x).
\end{equation}
Also, define the function $h_t:\mathcal{X}\to [0,1]$ by
\begin{equation}\label{htdef}
h_t(x):=\E[H_t(x)|\D],
\end{equation}
and then Fubini's theorem gives
\begin{equation}\label{intoverx}
\E[\Err_{t,1}\big| \D] = \ts\int_{\mathcal{X}} h_t(x) d\mu_1(x).
\end{equation}
A special property of the function $h_t(x)$ is that it only depends on $x$ through $\theta(x)$. To see this,  first let
$U_1,\dots,U_t$ be i.i.d.~Uniform[0,1] variables (independent of the objects $\D, \{Q_i\}$, and $(X,Y)$), and define the function $g_t:[0,1]\to [0,1]$ according to
\begin{equation}\label{gtdef}
g_t(\theta_0):=\P\Big(\ts\frac{1}{t}\tsum_{i=1}^t 1\{U_i\leq \theta_0\}\leq \ts\frac{1}{2}\Big),
\end{equation}
for any $\theta_0\in[0,1]$.
Since the sequences $\{Q_i(x)\}$ and  $\big\{1\{U_i\leq \theta(x)\}\big\}$ are both i.i.d.~Bernoulli$(\theta(x))$, conditionally on $\D$, it follows that we have the identity
\begin{equation}\label{gtidentity}
\begin{split}
h_t(x) =
g_t(\theta(x)),
\end{split}
\end{equation}
for all $x$ and $t$.
Consequently,  we may change variables from $x$ to $\theta(x)$ in line~\eqref{intoverx}, and then integrate over the unit interval to obtain
 \begin{equation*}
 \begin{split}
\E[\Err_{t,1}\big| \D] 
&=\int_0^1 g_t(\theta) f_1(\theta)d\theta\\
&=\int_0^{1/2} \Big(g_t(\theta)f_1(\theta)+g_t(1-\theta)f_1(1-\theta)\Big)d\theta,
\end{split}
\end{equation*}
where, in the second line, we have replaced $\theta$ with $(1-\theta)$ over the half interval $[1/2,1]$.
 To simplify things a bit further, note that because $F_1(\frac{1}{2})=\int_0^{1/2}f_1(\theta)d\theta$, the difference $\E[\Err_{t,1}\big| \D]-F_1(\frac{1}{2})$ can be written as
 \begin{equation*}
 \begin{split}
\E[\Err_{t,1}\big| \D] -F_1(\ts\frac{1}{2})
&=\int_0^{1/2} \Big(\big(g_t(\theta)-1\big)f_1(\theta)+g_t(1-\theta)f_1(1-\theta)\Big)d\theta.\\
\end{split}
\end{equation*}
Next, we use a special property of the function $g_t$. Specifically, if we let $G_{t,p}(\ts\frac{t}{2})$ denote the binomial c.d.f.~evaluated at $t/2$ (based on $t$ trials with success probability $p$), it is simple to check that the relation
$G_{t,1-p}(\ts\frac{t}{2})=1-G_{t,p}(\ts\frac{t}{2})$ holds for all $p\in[0,1]$ and odd $t$.
In terms of the function $g_t$, this means
\begin{equation}\label{gsym}
 g_t(1-\theta)=1-g_t(\theta),
\end{equation}
for all $\theta\in [0,1]$ and odd $t$. Consequently, the previous integral becomes
 \begin{equation*}
 \begin{split}
\E[\Err_{t,1}\big| \D] -F_1(\ts\frac{1}{2})
&=\int_0^{1/2} g_t(1-\theta)\Big(f_1(1-\theta)-f_1(\theta)\Big)d\theta.\\
\end{split}
\end{equation*}
 The quantity $f_1'(\frac{1}{2})$ now emerges by changing variables from $\theta$ to a new variable $u$ via the relation $\theta=\frac{1}{2}-\frac{u}{\sqrt{t}}$, with $u$ ranging over the interval $[0,\sqrt{t}/2]$. When we scale the difference $\E[\Err_{t,1}\big| \D] -F_1(\ts\frac{1}{2})$ by a factor of $t$, we obtain
 \begin{equation}\label{mainrep}
 \footnotesize
t\,\Big(\E[\Err_{t,1}|\D] -F_1(\ts\frac{1}{2})\Big)= \displaystyle \int_0^{\sqrt{t}/2}\! \underbrace{ g_t(\ts\frac{1}{2}+\ts\frac{u}{\sqrt{t}})\cdot \ts\frac{1}{2u/\sqrt{t}}\Big(\displaystyle f_1(\ts\frac{1}{2}+\ts\frac{u}{\sqrt{t}})-f_1(\ts\frac{1}{2}-\ts\frac{u}{\sqrt{t}})\Big)\cdot 2u\,}_{=:\, \psi_t(u)} du,
 \end{equation}
where we have defined the function $\psi_t$ above. Note also that a factor of $\sqrt{t}$ is absorbed by the relation $d\theta=-\ts\frac{1}{\sqrt{t}}du$.\\

To finish the proof, we evaluate the pointwise limit of $\psi_t$ and apply the dominated convergence theorem. (The task of showing that $\psi_t(u)$ is dominated by a suitable sequence of functions will be handled in a separate paragraph at the end of this subsection.) To proceed,  we first compute the pointwise limit of $\psi_t(u)$.
Since $f_1$ is differentiable at $1/2$, it is clear that as $t\to\infty$
 \begin{equation}\label{derivlimit}
 \ts\frac{f(\ts\frac{1}{2}+\ts\frac{u}{\sqrt{t}})-f(\ts\frac{1}{2}-\ts\frac{u}{\sqrt{t}})}{2u/\sqrt{t}}\ \to \  f_1'(\ts\frac{1}{2}).
 \end{equation}
  The limit of $g_t(\ts\frac{1}{2}+\ts\frac{u}{\sqrt{t}})$ is less obvious, and we compute it by expressing $g_t$ in terms of an empirical process. Letting the variables $U_1,\dots,U_t$ be as before, we define the random distribution function $\F_t(\theta) := \ts\frac{1}{t}\sum_{i=1}^t 1\{U_i \leq \theta\}$ for any $\theta\in [0,1]$. This gives
\begin{equation*}
\begin{split}
g_t(\ts\frac{1}{2}+\ts\frac{u}{\sqrt{t}}) &= \P\Big\{ \F_t(\ts\frac{1}{2}+\ts\frac{u}{\sqrt{t}})\leq \ts\frac{1}{2}\Big\}\\[0.2cm]
&= \P\Big\{ \sqrt{t}\Big(\F_t(\ts\frac{1}{2}+\ts\frac{u}{\sqrt{t}}) -(\ts\frac{1}{2}+\ts\frac{u}{\sqrt{t}}) \Big)  \leq -u\Big\},
\end{split}
\end{equation*}
where the second line involves a bit of algebra.
In order to evaluate the limit of this expression, we use a consequence of Donsker's Theorem~\cite[Lemma 19.24]{vaart}. Namely,  if $\{\theta_t\}\subset [0,1]$ is a numerical sequence that converges to a constant $\theta_0$, then the following limit in distribution holds
\begin{equation*}
\sqrt{t}\big(\F_t(\theta_t)-\theta_t\big)\xrightarrow{ \   \ d \  \ } N\big(0, \,\theta_0(1-\theta_0)\big).
\end{equation*}
Next, consider taking $\theta_t=\ts\frac{1}{2}-\ts\frac{u}{\sqrt{t}}$ and $\theta_0=1/2$. Also, let $Z\sim N(0,1)$ so that $\ts\frac{1}{2}Z\sim N(0,\ts\frac{1}{4})$. Therefore, with $u$ held fixed, the previous limit implies that as $t\to\infty$,
\begin{equation}\label{gtlimit}
g_t(\ts\frac{1}{2}+\ts\frac{u}{\sqrt{t}}) \xrightarrow{ \ \  \ }  \P(\ts\frac{1}{2}Z\leq -u)=\Phi(-2u),
\end{equation}
where $\Phi$ is the standard normal distribution function. Combining this with line~\eqref{derivlimit}, and the definition of $\psi_t$ in line~\eqref{mainrep}, we have the pointwise limit
$$\psi_t(u) \to  2\cdot f_1'(\ts\frac{1}{2})\cdot  u\cdot  \Phi(-2u).$$
So, provided that the dominated convergence theorem may be applied to $\psi_t(u)$, we conclude that
\begin{equation*}
\begin{split}
\lim_{t\to\infty} t\big(\E[\Err_{t,1}|\D]-F_1(\ts\frac{1}{2})\big) &= \displaystyle  2f_1'(\ts\frac{1}{2})\cdot  \int_0^{\infty}\, u\,\Phi(-2u)\, du\\[0.2cm]
&=\ts\frac{1}{8} f_1'(\ts\frac{1}{2}),
\end{split}
\end{equation*}
as needed (where the last line follows from a short integration-by-parts calculation).\qed \\
\paragraph{Details for showing $\psi_t(u)$ is dominated} We use a slightly generalized version of the standard dominated convergence theorem~\cite[Theorem 1.21]{kallenberg}. In particular, it is enough to construct a sequence of non-negative functions $b(u),b_1(u),b_2(u),\dots$ that are integrable on $[0,\infty)$, and satisfy the following three conditions
\begin{align}
 |\psi_t(u)|\cdot  1\{ u\leq \sqrt{t}/2\} & \  \leq \ b_t(u)\label{firstcondition}\\[0.3cm]
  b_t(u) & \ \to \ b(u)\label{secondcondition}\\[0.3cm]
 \ts\int_{0}^{\infty}b_t(u)du & \ \to \ \ts\int_0^{\infty}b(u)du.\label{thirdcondition}
\end{align}
 The main idea is now to bound $\psi_t(u)$ in two pieces, depending on the size of $u$. Due to \textbf{A2}, there are constants $\kappa>0$ and $\delta_1\in(0,1/2)$ such that the difference quotient of $f_1$ satisfies the following bound for every $t\geq 1$,
\begin{equation*}
\bigg|\ts\frac{\displaystyle f_1(\ts\frac{1}{2}+\ts\frac{u}{\sqrt{t}})-f_1(\ts\frac{1}{2}-\ts\frac{u}{\sqrt{t}})}{2u/\sqrt{t}} \bigg|\
\leq \kappa  \ \ \ \ \text{ when  } \ \ \ \ u \leq \delta_1\sqrt{t}.
\end{equation*}
Next, in order to control $g_t(\frac{1}{2}+\frac{u}{\sqrt{t}})$, we apply Hoeffding's inequality~\cite[Theorem 2.8] {boucheron2013} to the binomial distribution and the definition of $g_t$ in line~\eqref{gtdef} to obtain
\begin{equation}\label{hoeffding1}
g_t(\ts\frac{1}{2}+\ts\frac{u}{\sqrt{t}})\leq e^{-2u^2},
\end{equation}
for all $0\leq u\leq \sqrt{t}/2$ and every $t\geq 1$. Likewise, we define the non-negative function
\footnotesize
\begin{equation*}
b_t(u):=\begin{cases}  & 2\kappa\cdot u\cdot  \exp(-2u^2) \ \ \ \ \ \ \ \ \ \ \ \ \ \   \ \ \ \ \  \ \, \ \ \  \ \ \ \ \ \  \  \ \ \ \  \ \ \ \ \text{ when }  \ \ \ \ u\in [0,\delta_1\sqrt{t}]\\
& \sqrt{t} \exp(-2\delta_1^2 t) \cdot \Big( f_1(\ts\frac{1}{2}+\ts\frac{u}{\sqrt{t}})+f_1(\ts\frac{1}{2}-\ts\frac{u}{\sqrt{t}})\Big)  \ \ \ \ \ \text{ when } \ \ \  \ u\in (\delta_1\sqrt{t},\sqrt{t}/2]\\[0.2cm]
& 0  \  \ \ \ \ \ \  \ \ \ \ \ \    \ \ \ \  \ \ \ \    \ \ \ \  \ \ \ \   \ \ \ \  \ \ \ \    \ \ \ \  \ \ \ \  \ \ \ \ \,   \ \  \ \ \ \ \ \ \ \ \  \text{ when } \ \ \  \ u\in (\sqrt{t}/2,\infty).
\end{cases}
\end{equation*}
\normalsize
and the last few steps give
$$|\psi_t(u)|\cdot 1\{u\leq \sqrt{t}/2\}\leq b_t(u).$$
It is also simple to check that if we define
$$b(u):=2\kappa \cdot u \cdot \exp(-2u^2),$$
then we have the pointwise limit,
$$b_t(u) \to b(u).$$
Finally, to check the third condition~\eqref{thirdcondition}, observe that the change of variable $\theta=\ts\frac{1}{2}-\ts\frac{u}{\sqrt{t}}$ gives
$$\ \int_{\delta_1\sqrt{t}}^{\sqrt{t}/2}\Big(f_1(\ts\frac{1}{2}+\ts\frac{u}{\sqrt{t}})+f_1(\ts\frac{1}{2}-\ts\frac{u}{\sqrt{t}})\Big)du = \sqrt{t}\displaystyle \int_0^{1/2-\delta_1}(f_1(1-\theta)+f_1(\theta))d\theta=\mathcal{O}(\sqrt{t}).$$
Consequently, integral of $b_t(u)$ on the middle interval $(\delta_1\sqrt{t},\sqrt{t}/2]$ does not matter asymptotically, since it is driven to 0 by the factor $\sqrt{t}\exp(-2\delta_1^2t)$. This implies that the condition~\eqref{thirdcondition} holds.

\subsubsection{Proof of Lemma 1 in the non-differentiable case}\label{nondiffcase}
When $f_1$ is not differentiable at 1/2, the proof may be repeated in the same way up to line~\eqref{mainrep}. At this stage, we then replace $\psi_t(u)$ with the bounding function $b_t(u)$ to obtain the inequality
\begin{equation*}
\small
t\,\big|\E[\Err_{t,1}|\D] -F_1(\ts\frac{1}{2})\big| \ \leq \  \ts\int_0^{\infty}b_t(u)du,
\end{equation*}
for every $t\geq 1$. (Note that the properties~\eqref{firstcondition},~\eqref{secondcondition}, and~\eqref{thirdcondition} of the function $b_t(u)$ only depend on $f_1$ being Lipschitz in a neighborhood of 1/2.)
Next, due to the condition~\eqref{thirdcondition}, the sequence of numbers $\ts\int_0^{\infty}b_t(u)du$ is bounded by some positive constant $C$, and hence
\begin{equation*}
\small
t\,\big|\E[\Err_{t,1}|\D] -F_1(\ts\frac{1}{2})\big| \ \leq \  C,
\end{equation*}
as desired.\qed

\subsection{Proof of Theorem~\ref{VARTHM} (variance bound)}\label{VARTHMproof}

Instead of proving the bound~\eqref{varbound} directly, it will be more convenient to bound a related quantity. If we think of the random variable $\Err_{t,1}$ as an ``estimator'' of the parameter $F_1(\ts\frac{1}{2})$, then the standard decomposition MSE = variance + $\text{bias}^2$ gives the relation
\begin{equation*}
\E\bigg[ \Big(\Err_{t,1}-F_1(\ts\frac{1}{2})\Big)^2\Big| \D\bigg]= \var(\Err_{t,1}\big| \D)+ \Big(\E[\Err_{t,1}\big| \D]-F_1(\ts\frac{1}{2})\Big)^2.
\end{equation*}
Next, if we multiply through by $t$, and use Lemma~\ref{EXPECLEMMA} to note that $t\big(\E[\Err_{t,1}\big| \D]-F_1(\ts\frac{1}{2})\big)^2=\mathcal{O}(\ts\frac{1}{t})$, then we conclude
\begin{equation*}
\underbrace{t\, \E\bigg[ \Big(\Err_{t,1}-F_1(\ts\frac{1}{2})\Big)^2\Big| \D\bigg]}_{(*)} = t\,\var(\Err_{t,1}\big| \D)+o(1).
\end{equation*}
Here, the $(*)$ symbol is merely a shorthand that will be convenient in the remainder of the proof. Thus, it is enough prove $(*)\leq \frac{1}{4}f_1(\ts\frac{1}{2})^2+o(1)$.

To begin with the main portion of the proof, define the complementary sets
$$\mathcal{X}_+:=\Big\{x\in \mathcal{X} : \, \theta(x)\leq \ts\frac{1}{2}\Big\} \ \ \text{ and } \ \ \mathcal{X}_-:=\Big\{x\in \mathcal{X} :\, \theta(x)> \ts\frac{1}{2}\Big\}. $$
Recalling the notation $\mu_1=\mathcal{L}(X|Y=1)$, note that by \textbf{A2} and the standard change-of-variable rule, we have
\begin{equation*}
\begin{split}
F_1(\ts\frac{1}{2})&=\ts\int_0^{1}1\{\theta\leq \ts\frac{1}{2}\} f_1(\theta)d\theta \\[0.3cm]
&=\ts\int_{\mathcal{X}} 1\{\theta(x)\leq \ts\frac{1}{2}\} d\mu_1(x)\\[0.3cm]
&=\ts\int_{\mathcal{X}_+} d\mu_1(x).
\end{split}
\end{equation*}
Combining this with the representation of $\Err_{t,l}$ in line~\eqref{Hrep}, it follows that
\begin{align}
(*) \ 
&= t\, \E\bigg[\Big(\ts\int_{\mathcal{X}_+}(H_t(x)-1)d\mu_1(x) +\ts\int_{\mathcal{X}_-} H_t(x)d\mu_1(x)\Big)^2\bigg|\D\bigg]\\[0.3cm]
&\leq t\, \E\bigg[ \Big(\ts\int_{\mathcal{X}_+} (H_t(x)-1)d\mu_1(x)\Big)^2\bigg| \D\bigg] +t\, \E\bigg[ \Big(\ts\int_{\mathcal{X}_-} H_t(x)d\mu_1(x)\Big)^2\bigg|\D\bigg],
\end{align}
where the inequality comes from dropping the cross-term, since the integral over $\mathcal{X}_+$ is at most 0, and the integral over $\mathcal{X}_-$ is at least 0. Next, we write the squared integrals as double integrals and use Fubini's theorem to obtain

\begin{equation*}
\begin{split}
(*) \ 
&\leq t\, \int_{\mathcal{X}_+}\int_{\mathcal{X}_+} \E\Big[(H_t(x)-1)(H_t(x')-1)\Big|\D\Big]d\mu_1(x)d\mu_1(x')\\[0.2cm]
& \ \ \ \  + t\, \int_{\mathcal{X}_-}\int_{\mathcal{X}_-} \E\Big[H_t(x)H_t(x')\Big|\D\Big]d\mu_1(x)d\mu_1(x').
\end{split}
\end{equation*}
Recall the function $h_t(x)=\E[H_t(x)|\D]$ from line~\eqref{htdef}. Due to the fact that $H_t(\cdot)$ is binary, we clearly have that the product $H_t(x)H_t(x')$ is upper-bounded by $H_t(x)$ and $H_t(x')$. It follows that the two integrands in the previous line can be bounded using 
\begin{equation}\label{hmin}
\E\Big[H_t(x)H_t(x')\Big|\D\Big] \leq \min\big\{h_t(x),h_t(x')\big\},
\end{equation}
which leads to the following inequality after expanding the product $(H_t(x)-1)(H_t(x')-1)$,
\begin{equation}\label{xminusint}
\begin{split}
(*) \ 
&\leq t\, \int_{\mathcal{X}_+}\int_{\mathcal{X}_+} \bigg(\min\big\{h_t(x),h_t(x')\big\}-h_t(x)-h_t(x')+1\bigg)d\mu_1(x)d\mu_1(x')\\[0.2cm]
& \ \ \ \ + t\, \int_{\mathcal{X}_-}\int_{\mathcal{X}_-} \min\big\{h_t(x),h_t(x')\big\}d\mu_1(x)d\mu_1(x').
\end{split}
\end{equation}
At this point, we make use of the identity $h_t(x)=g_t(\theta(x))$, derived in line~\eqref{gtidentity} in the proof of Lemma~\ref{EXPECLEMMA}. Due to \textbf{A2}, we may use a change of variable to integrate the density $f_1$ over the unit interval, rather than integrating $\mu_1$ over $\mathcal{X}$.
In particular, note that the sets $\mathcal{X}_+$ and $\mathcal{X}_-$ correspond to the half intervals $[0,\frac{1}{2}]$ and $(\frac{1}{2},1]$ respectively. Furthermore, for the second integral in line~\eqref{xminusint},  we can make another change of variable by replacing $\theta$ with $1-\theta$, so that all integrals are over $[0,\frac{1}{2}]$, which leads to
\small
\begin{align*}
(*) \ 
&\leq t\, \int_0^{1/2}\int_0^{1/2}\Bigg( \min\big\{g_t(\theta),g_t(\theta')\big\}-g_t(\theta)-g_t(\theta')+1\Bigg) f_1(\theta)f_1(\theta')d\theta\,d\theta'\\[0.2cm]
& \ \ \ \  + t\, \int_{0}^{1/2}\int_{0}^{1/2} \min\big\{g_t(1-\theta),g_t(1-\theta')\big\}f_1(1-\theta)f_1(1-\theta')d\theta\,d\theta'.
\end{align*}
\normalsize
It turns out that quite a bit of additional simplification is possible. First, note the simple identity%
$$ g_t(\theta)+g_t(\theta')= \min\big\{g_t(\theta),g_t(\theta')\big\}+\max\big\{g_t(\theta),g_t(\theta')\big\}.$$
Next, we use the fact that $g_t(1-s)=1-g_t(s)$ for all $s\in[0,1]$ (as was shown in line~\eqref{gsym}) to conclude
\begin{equation*}
\begin{split}
\min\big\{g_t(\theta),g_t(\theta')\big\}-g_t(\theta)-g_t(\theta')+1 \ & =   1-\max\{g_t(\theta),g_t(\theta')\}\\[0.2cm]
&=\min\{g_t(1-\theta),g_t(1-\theta')\}.
\end{split}
\end{equation*}
Hence, the previous integrals can be combined as
\footnotesize
$$
(*) \ 
\leq t\, \int_0^{1/2}\int_0^{1/2}\bigg( \min\big\{g_t(1-\theta),g_t(1-\theta')\}\bigg) \bigg(f_1(\theta)f_1(\theta')+f_1(1-\theta)f_1(1-\theta')\bigg)d\theta\,d\theta'.
$$
\normalsize
Now consider the change of variable  $(\theta,\theta')=(\ts\frac{1}{2}-\frac{u}{\sqrt{t}},\ts\frac{1}{2}-\frac{u'}{\sqrt{t}})$ with $u$ and $u'$ ranging over $[0,\sqrt{t}/2]$, and note that a factor of $t$ is absorbed by the relation $t \, d\theta\,d\theta' = du\,du'$. Likewise, defining the function
\begin{equation}\label{gammatdef}
\footnotesize
\gamma_t(u,u')  \ := \ \min\Big\{g_t(\ts\frac{1}{2}+\frac{u}{\sqrt{t}}),g_t(\ts\frac{1}{2}+\frac{u'}{\sqrt{t}})\Big\}\bigg(f_1(\ts\frac{1}{2}-\frac{u}{\sqrt{t}})f_1(\ts\frac{1}{2}-\frac{u'}{\sqrt{t}})+ f_1(\ts\frac{1}{2}+\frac{u}{\sqrt{t}})f_1(\ts\frac{1}{2}+\frac{u'}{\sqrt{t}})\bigg)
\end{equation}
gives
$$
(*)\leq  \int_{0}^{\sqrt{t}/2}\int_{0}^{\sqrt{t}/2}\gamma_t(u,u') du\,du'.
$$
Using the limit~\eqref{gtlimit} from the proof of Lemma~\ref{EXPECLEMMA}, and the continuity of $f_1$ at 1/2, it follows that if we define
\begin{equation*}
\begin{split}
\gamma(u,u') \ &:= 2f_1(\ts\frac{1}{2})^2\cdot \min\big\{\Phi(-2u),\Phi(-2u')\big\}\\[0.2cm]
& \  = 2f_1(\ts\frac{1}{2})^2\cdot\Phi(-2\max\{u,u'\}) ,
\end{split}
\end{equation*}
then we have the pointwise limit
$$\gamma_t(u,u') \ \to \  \gamma(u,u').$$
So, provided that this limit is dominated (which will be handled at the end of this subsection), the dominated convergence theorem yields

\begin{equation}\label{gammaint}
(*)\leq \int_0^{\infty}\int_0^{\infty} \gamma(u,u')dudu' + o(1).
\end{equation}
To compute this integral, let $\mathcal{R}_+$ denote the set of pairs $(u,u')$ in the quadrant $[0,\infty)^2$ such that $u'\geq u$, and let $\mathcal{R}_-$ denote the set of pairs where $u'< u$. Then,
\begin{equation*}
\footnotesize
\begin{split}
\int_0^{\infty}\int_0^{\infty} \gamma(u,u')dudu' 
&=2f_1(\ts\frac{1}{2})^2 \displaystyle\int_{\mathcal{R}_+}\Phi(-2u') du du' + 2f_1(\ts\frac{1}{2})^2\displaystyle \int_{\mathcal{R}_-}\Phi(-2u)du du' \\[0.3cm]
&=4f_1(\ts\frac{1}{2})^2 \displaystyle\int_{\mathcal{R}_+}\Phi(-2u') du du'\ \ \ \ \text{(by symmetry)}\\[0.3cm]
&=4f_1(\ts\frac{1}{2})^2 \displaystyle\int_0^{\infty} u'\Phi(-2u') du' \\[0.3cm]
&=\ts\frac{1}{4}f_1(\ts\frac{1}{2})^2,
\end{split}
\end{equation*}
as desired, where the last line follows from an integration-by-parts calculation.\qed

\paragraph{Details for showing $\gamma_t(u,u')$ is dominated} 
To apply the dominated convergence theorem 
~\cite[Theorem 1.21]{kallenberg}, it is enough to construct non-negative functions $b(u,u'),b_1(u,u'),b_2(u,u'),\dots$ that are integrable on the quadrant $[0,\infty)^2$, and satisfy the following three conditions,
\begin{align}
| \gamma_t(u,u') | \cdot 1\big\{u\leq \sqrt{t}/2,u'\leq \sqrt{t}/2\big\}& \  \leq \ b_t(u,u')\label{firstconditionagain}\\[0.3cm]
  b_t(u,u') & \ \to \ b(u,u')\label{sconditionagain}\\[0.3cm]
 \ts\int_0^{\infty}\int_{0}^{\infty}b_t(u,u')dudu' & \ \to \ \ts\int_0^{\infty}\int_0^{\infty}b(u,u')dudu'.\label{thirdconditionagain}
\end{align}
To construct these functions, recall our assumption that $f_1$ is Lipschitz in a neighborhood of $1/2$, and so there are positive constants $\delta_1$ and $\kappa$ such that
$$f_1(\ts\frac{1}{2}\pm \ts\frac{u}{\sqrt{t}}) \leq f(1/2)+\kappa\delta_1 \ \ \ \ \ \ \ \text{ when }  \ \ \ \ \ u\leq \delta_1\sqrt{t}.$$
Also, recall the Hoeffding bound from line~\eqref{hoeffding1},
\begin{equation*}
g_t(\ts\frac{1}{2}+\ts\frac{u}{\sqrt{t}})\leq e^{-2u^2},
\end{equation*}
which holds for all $0\leq u\leq \sqrt{t}/2$ and every $t\geq 1$.
Accordingly, by looking the definition of $\gamma_t(u,u')$ in line~\eqref{gammatdef} we define the bounding function
%
\scriptsize
\begin{equation*}
b_t(u,u') :=
\begin{cases}
& \!\!\!\!2\, e^{-2\max\{u^2,u'^2\}}  \Big(f_1(\ts\frac{1}{2})+\delta_1\kappa\Big)^2  \ \ \ \ \ \ \ \ \ \ \ \ \ \   \ \ \ \ \ \ \  \ \  \  \text{ when } \ \ \ \  \ u,u'\leq \delta_1\sqrt{t}\\[0.3cm]
&  \!\!\!\!e^{-2 \delta_1^2 t}\Big(f_1(\ts\frac{1}{2})+\delta_1\kappa\Big)\Big( f_1(\ts\frac{1}{2}-\frac{u}{\sqrt{t}})+f_1(\ts\frac{1}{2}+\frac{u}{\sqrt{t}})\Big)  \ \    \text{ when }    \ u'\leq \delta_1\sqrt{t} \ \ \text{ and } u\in(\delta_1\sqrt{t},\sqrt{t}/2]\\[0.3cm]
& \!\!\!\! e^{-2 \delta_1^2 t}\Big(f_1(\ts\frac{1}{2})+\delta_1\kappa\Big)\Big( f_1(\ts\frac{1}{2}-\frac{u'}{\sqrt{t}})+f_1(\ts\frac{1}{2}+\frac{u'}{\sqrt{t}})\Big)  \  \  \text{ when } \ u\leq \delta_1\sqrt{t} \ \ \text{ and } u'\in(\delta_1\sqrt{t},\sqrt{t}/2]\\[0.3cm]
&\!\!\!\! 0  \ \  \  \  \    \ \  \  \  \       \  \ \  \   \    \  \   \  \ \ \  \  \  \    \  \   \  \ \  \ \  \  \  \    \  \   \  \ \   \  \   \  \ \  \ \  \  \  \    \  \   \  \ \  \ \  \  \  \  \ \ \ \ \ \ \ \text{ when }\ \ \  \  \ u, u'>\delta_1\sqrt{t}.
\end{cases}
\end{equation*}
\normalsize
It is straightforward to check the bounding condition~\eqref{firstconditionagain}. Furthermore, if we define
\begin{equation*}
b(u,u') :=
 2\, e^{-2\max\{u^2,u'^2\}}  \big(f_1(\ts\frac{1}{2})+\delta_1\kappa \big)^2,
\end{equation*}
then we have the following pointwise limit for fixed $u$ and $u'$,
$$b_t(u,u') \to b(u,u').$$
Finally, we check the third condition~\eqref{thirdconditionagain} for dominated convergence. The essential point to notice is that the second and third lines in the definition of $b_t$ are asymptotically negligible when integrating over $u$ and $u'$. To see this, consider the second line in the definition of $b_t$, and note that the change of variable $\theta=\ts\frac{1}{2}-\ts\frac{u}{\sqrt{t}}$  implies
\begin{equation*}
\begin{split}
2\, e^{-2\delta_1^2t}\int_0^{\delta_1\sqrt{t}}\int_{\delta_1\sqrt{t}}^{\sqrt{t}/2}f_1(\ts\frac{1}{2}-\frac{u}{\sqrt{t}})dudu'   \ & =  \ \delta_1 t  \, e^{-2\delta_1^2t}\int_0^{1/2-\delta_1}f_1(\theta)d\theta\\[0.2cm]
&=o(1).
\end{split}
\end{equation*}
\normalsize
Likewise, by repeating this calculation using the function $f_1(\ts\frac{1}{2}+\frac{u}{\sqrt{t}})$, as well as by interchanging the roles of $u$ and $u'$,  the condition~\eqref{thirdconditionagain} follows.

\subsection{Proof of Theorem~\ref{ATTAINTHM} (attaining the variance bound)}\label{ATTAINTHMproof}

 Note that $\bar{Q}^{\circ}_t(x) = \F_t(\theta(x))$, where again $\F_t(\cdot) = \frac{1}{t}\sum_{i=1}^t 1\{U_i\leq \cdot\}$, and the variables $U_1,\dots,U_t$ are the same as in the definition of $\{Q_i^{\circ}(\cdot)\}$.  Due to \textbf{A2}, we may use a change of variable to express $\Err_{t,1}^{\circ}$ as
\begin{equation*}
\Err^{\circ}_{t,1} = \int_0^1 1\{\F_t(\theta)\leq \ts\frac{1}{2}\} f_1(\theta)d\theta.
\end{equation*}
For any $r\in(0,1)$,  define the quantile function
$\F_t^{-1}(r) := \inf\{\theta: \F(\theta)\geq r\}$, and recall the equivalence $\F_t(\theta)< r \Longleftrightarrow\theta<\F_t^{-1}(r) $ for any $\theta\in[0,1]$ (see~\cite[Lemma 21.1]{vaart}). With this fact in hand, the variable $\Err_{t,1}^{\circ}$ can be evaluated as follows, where we note that $1\{\F_t(\theta)\leq \ts\frac{1}{2}\}=1\{\F_t(\theta)< \ts\frac{1}{2}\}$ when $t$ is odd,
\begin{equation*}
\begin{split}
\Err^{\circ}_{t,1} 
%
&=\int_0^1 1\{\F_t(\theta)<\ts\frac{1}{2}\} f_1(\theta)d\theta \text{  }\\[0.2cm]
&=\int_0^1 1\{\theta < \ts\F_t^{-1}(\frac{1}{2})\} f_1(\theta)d\theta\\[0.2cm]
&=\int_0^{\F_t^{-1}(\frac{1}{2})}  f_1(\theta)d\theta\\[0.2cm]
&=F_1(\F_t^{-1}(\ts\frac{1}{2})).
\end{split}
\end{equation*}
Next, we use the fact that the quantile process $\F_t^{-1}(\frac{1}{2})$ satisfies the following  limit in distribution~\cite[Corollary 21.5]{vaart},
\begin{equation*}
\sqrt{t}(\F_t^{-1}(\ts\frac{1}{2})-\ts\frac{1}{2}) \xrightarrow{ \ \ d \ \ } N(0,\frac{1}{4}).
\end{equation*}
Also, when $f_1$ is Lipschitz in a neighborhood of 1/2, it follows that $F_1$ is differentiable at $1/2$, and then the delta method~\cite[Theorem 3.1]{vaart} gives
\begin{equation}\label{delta}
\sqrt{t}\Big(\Err^{\circ}_{t,1}-F_1(\ts\frac{1}{2})\Big)   \ = \sqrt{t}\Big(F_1(\F_t^{-1}(\ts\frac{1}{2}))-F_1(\ts\frac{1}{2})\Big)\xrightarrow{ \ \ d \ \ } N\Big(0, \, \ts \ts\frac{1}{4}f_1(\ts\frac{1}{2})^2 \Big).
\end{equation}
From Lemma~\ref{EXPECLEMMA}, we know that 
$$\sqrt{t}\big(\E[\Err^{\circ}_{t,l}\big| \D]-F_1(\ts\frac{1}{2}))=o(1),$$
 and if we define the zero-mean random variable
$$\zeta_t:=\sqrt{t}\big(\Err^{\circ}_{t,l}-\E[\Err^{\circ}_{t,l}\big| \D]\big),$$
 then
it follows from Slutsky's lemma~\cite[Lemma 2.8]{vaart} that $\zeta_t$ satisfies the same distributional limit as in line~\eqref{delta}, namely
$$\zeta_t\xrightarrow{ \ \ d \ \ } N\Big(0,\ts\frac{1}{4}f_1(\ts\frac{1}{2})^2\Big).$$
Finally, it is a general fact that if a sequence of zero-mean random variables has a distributional limit $\zeta_t\xrightarrow{ \ d \ } \zeta$, then the limiting variance satisfies $\var(\zeta)\leq \liminf_{t\to\infty} \var(\zeta_t)$~\cite[Lemma 4.11]{kallenberg}. Hence,
\begin{equation*}
\ts\frac{1}{4}f_1(\ts\frac{1}{2})^2 \ \leq \ \displaystyle \liminf_{t\to\infty} \, t\,\var(\Err^{\circ}_{t,1}|\D).
\end{equation*}
On the other hand, the upper bound in Theorem~\ref{VARTHM} gives 
\begin{equation*}
 \limsup_{t\to\infty}\,  t\,\var(\Err^{\circ}_{t,1}|\D) \ \leq  \ \ts\frac{1}{4}f_1(\ts\frac{1}{2})^2,
\end{equation*}
and so combining the last two statements leads to the desired limit.\qed

\subsection{Proof of Theorem~\ref{MSETHM} (MSE bound)}\label{MSETHMproof}

 For each $j=1,\dots,m_l$, define the random variable
\begin{equation*}
\Delta_{j}:=\ts\frac{1}{h}K\Big(\ts\frac{\bar{Q}_t(\tilde{X}_{j,l})-1/2}{h}\Big)-f_l(\ts\frac{1}{2}).
\end{equation*}
(To ease notation, we suppress the fact that $\Delta_j$ depends on $h,t,l,$ and $m_l$.)
Then,
\begin{equation*}
\hat{f}_l(\ts\frac{1}{2})-f_l(\ts\frac{1}{2}) = \ts\frac{1}{m_l}\displaystyle\sum_{j=1}^{m_l} \Delta_j, 
\end{equation*}
and since the hold-out points $\tilde{X}_{1,l},\dots,\tilde{X}_{m_l,l}$ are i.i.d., we have
\begin{equation*}
\textsc{mse}\big(\hat{f}_l(\ts\frac{1}{2})\big) = \ts\frac{1}{m_l} \E[\Delta_{1}^2|\D]+ \ts\frac{2}{m_l^2}\binom{m_l}{2} \E[\Delta_{1}\Delta_{2}\big| \D].
\end{equation*}
The remainder of the proof deals with the task of deriving bounds for $\E[\Delta_{1}^2\big| \D]$ and $\E[\Delta_{1}\Delta_{2}\big| \D]$, as addressed below in Lemmas~\ref{squarelemma} and~\ref{crosslemma} (respectively).
The theorem then follows by choosing the bandwidth $h$  that minimizes the sum of the bounds (in terms of rates), as described in~\ref{bandwidthexpl} below.
\qed

\begin{lemma}\label{squarelemma}
Suppose the conditions of Theorem~\ref{MSETHM} hold, and let $[1/2-\delta_l,1/2+\delta_l]$ denote an interval on which $f_l$ or $f_l'$ is Lipschitz, with $\delta_l\in(0,1/2)$. Then, there is a number $\kappa>0$ not depending on $t$ or $m_l$, such that for any $h\in(0,\delta_l)$,
\begin{equation*}
\ts\frac{1}{m_l}\E[\Delta_1^2\big| \D]\leq \ts\frac{\kappa}{h \, m_l}\Big(1+\ts\frac{1}{h\sqrt{t}}\Big).
\end{equation*}
\end{lemma}

\paragraph{Remark} To simplify the statement of the following lemma, we will refer to $f_l$ and $f_l'$ as $f_l^{(\beta-1)}$ with $ \beta=1$ or $\beta=2$ (respectively).
\begin{lemma}\label{crosslemma}
Fix $\beta\in\{1,2\}$. Suppose the conditions of Theorem~\ref{MSETHM} hold, and let $[1/2-\delta_l,1/2+\delta_l]$ denote an interval on which $f_l^{(\beta-1)}$ is Lipschitz, with $\delta_l\in(0,1/2)$. Then, there is a number $\kappa>0$, not depending on $t$ or $m_l$ such that for any $h\in(0,\delta_l)$,
\begin{equation*}
\ts\frac{2}{m_l^2}\binom{m_l}{2} \E[\Delta_1\Delta_2|\D]\leq \kappa\big(h^{2\beta}+\ts\frac{1}{h\sqrt{t}}+\ts\frac{1}{h^2t}\big).
\end{equation*}
\end{lemma}

\paragraph{Remark} The proof of Lemmas~\ref{squarelemma} and~\ref{crosslemma} are given in~\ref{sec:squarelemma} and~\ref{sec:crosslemma} below. The proofs will only be given for the case $l=1$, since the proofs for $l=0$ are essentially the same.

\subsubsection{Explanation of bandwidth choice}\label{bandwidthexpl}
To explain our choice of bandwidth, we aim to express $h$ as a function of $t$ and $m_l$ so that the sum 
$$\ts\frac{1}{hm_l}\Big(1+\ts\frac{1}{h\sqrt{t}}\Big)+\Big(h^{2\beta}+\ts\frac{1}{h\sqrt{t}}+\ts\frac{1}{h^2t}\Big)$$
decreases at the fastest possible rate as $t,m_l\to\infty$ and  $h\to 0$.
As a first observation, note that the term $\ts\frac{1}{h^2t}$ can be dropped, because it will always be of smaller order than $\ts\frac{1}{h\sqrt{t}}$, provided that the latter quantity tends to 0. The same reasoning allows the term $\frac{1}{h m_l}\frac{1}{h\sqrt{t}}$ to be dropped as well.
Hence, it is enough to optimize the rate for the quantity
$$h^{2\beta}+\ts\frac{1}{hm_l}+\ts\frac{1}{h\sqrt{t}}.$$
Another simplification can be made by noting that the sum $\ts\frac{1}{hm_l}+\ts\frac{1}{h\sqrt{t}}$ will have the same rate as the slower of the two terms, which has the same rate as $1/(h \min\{\sqrt{t},m_l\})$. Finally, since the quantity $1/(h \min\{m_l,\sqrt{t}\})$ increases as $h$ becomes small, and the quantity $h^{2\beta}$ decreases as $h$ becomes small, the best choice of $h$ occurs when both quantities have matching rates. This leads to solving the rate equation
$$h^{2\beta} \asymp \ts\frac{1}{h}\cdot \frac{1}{\min\{\sqrt{t},m_l\}},$$
yielding
$$h=c_1\big(\min\{m_l,\sqrt{t}\}\big)^{\frac{-1}{2\beta+1}},$$
for some constant $c_1>0$, which is the stated bandwidth choice in Theorem 3.

\subsubsection{Proof of Lemma~\ref{squarelemma}}\label{sec:squarelemma}
  For the rectangular kernel $K(\cdot)=\ts\frac{1}{2}1\{-1\leq \cdot \leq 1\}$, we have the relation $(\frac{1}{h}K(\cdot))^2= \frac{1}{2h^2}K(\cdot)$, and some arithmetic leads to
\begin{equation*}
\begin{split}
\ts\frac{1}{m_1}\E[\Delta_1^2|\D] &= \ts\frac{1}{m_1h} \E\Big[\big(\ts\frac{1}{2h}-2f_1(\ts\frac{1}{2})\big)\cdot K\Big(\ts\frac{\bar{Q}_t(\tilde{X}_{1,1})-1/2}{h}\Big)\Big| \D\Big]+\frac{1}{m_1}f_1(\ts\frac{1}{2})^2\\[0.4cm]
& \leq  \ts\frac{1}{2m_1h} \E\Big[\ts\frac{1}{h} K\Big(\ts\frac{\bar{Q}_t(\tilde{X}_{1,1})-1/2}{h}\Big)\Big|\D\Big]+\frac{1}{m_1}f_1(\ts\frac{1}{2})^2,
\end{split}
\end{equation*}
where the inequality comes from dropping the negative term involving $-2f_1(\ts\frac{1}{2})$.

In the rest of the proof, it is enough to show there is a constant $\kappa>0$ not depending on $t$ or $m_1$, such that for any $h\in (0,\delta_1)$, 
 \begin{equation}\label{kappabound}
 \E\Big[\ts\frac{1}{h} K\Big(\ts\frac{\bar{Q}_t(\tilde{X}_{1,1})-1/2}{h}\Big)\Big|\D\Big]\leq \kappa\Big(1+\ts\frac{1}{h\sqrt{t}}\Big).
 \end{equation}
To proceed, let $U_1,\dots,U_t$ be i.i.d.~Uniform[0,1] random variables, and for any numbers $s\in[0,1]$, and  $h\in(0,\delta_1)$, define the function
\begin{equation}\label{varphiprob}
\varphi_t(s;h):= \P\Big(\ts\frac{1}{2}-h\leq \ts\frac{1}{t}\tsum_{i=1}^t 1\{U_i\leq s\}\leq \ts\frac{1}{2}+h\Big).
\end{equation}
It is simple check the relation
\begin{equation}\label{Kexpec}
\ts\frac{1}{2}\varphi_t(\theta(x);h)=\E\Big[K\Big(\ts\frac{\bar{Q}_{t}(x)-1/2}{h}\Big)\Big|\D\Big],
\end{equation}
for all $x\in\mathcal{X}$, and all $h>0$.
Using \textbf{A2}, and a change of variable from $x$ to $\theta$,  we may take the expectation over $\tilde X_{1,1}$ as
\begin{equation*}
\E\Big[\ts\frac{1}{h} K\Big(\ts\frac{\bar{Q}_t(\tilde{X}_{1,1})-1/2}{h}\Big)\Big|\D\Big] = \ts\frac{1}{2h} \displaystyle\int_0^1\varphi_t(\theta;h)f_1(\theta)d\theta.
\end{equation*}
To handle the integral on the right side, define the following upper-bounding function, \mbox{$\mathscr{U}_t: [0,1]\to \R$,} which acts as an approximate indicator on $[1/2-h,1/2+h]$,
\begin{equation}\label{udef}
 \mathscr{U}_t(\theta;h):=
\begin{cases}
& e^{-2t\,(1/2-h-\theta)^2}, \text{  \ \ \ \ \ \ \ \ \ \  \ \ \ \ for } \theta< \ts\frac{1}{2}-h\\
&1,  \text{ \ \ \ \ \ \ \ \ \ \ \ \ \ \ \ \   \ \ \ \ \ \ \ \ \ \ \  \ \  \ \ for } \theta\in [\ts\frac{1}{2}-h,\ts\frac{1}{2}+h]\\
&  e^{-2t\,(\theta- (1/2+h))^2}, \text{ \ \ \ \ \ \! \ \ \ \ \ \ \  for } \theta> \ts\frac{1}{2}+h.
\end{cases}
\end{equation}
Furthermore, if we apply Hoeffding's inequality to the probability in line~\eqref{varphiprob}, it follows that
\begin{equation}\label{uupper}
\varphi_t(\theta;h)\leq \mathscr{U}_t(\theta;h),
\end{equation}
for all $\theta\in [0,1]$, and all $h\in(0,\delta_1)$. Integrating this bound over $[0,1]$ gives 
\footnotesize
\begin{equation}\label{threepieces}
\begin{split}
\E\Big[\ts\frac{1}{h} K\Big(\ts\frac{\bar{Q}_t(\tilde{X}_{1,1})-1/2}{h}\Big)\Big] &\leq %
\ts\frac{1 }{2 h}\displaystyle \int_0^{1/2-h} e^{-2t\,(1/2-h-\theta)^2}f_1(\theta)d\theta\\[0.2cm]
& \  \ \ \ \ \ \ + \ts\frac{F_1(1/2+h)-F_1(1/2-h)}{2h}\\[0.2cm]
&  \  \ \ \ \ \ \ + \ts\frac{1 }{2 h}\displaystyle \int_{1/2+h}^{1}   e^{-2t\,(1/2+h-\theta)^2}f_1(\theta) d\theta.
\end{split}
\end{equation}
\normalsize
For future reference, we define the middle term as the ``central difference quotient''
\begin{equation}\label{cdqdef}
 \textsc{cdq}(h) :=\ts\frac{F_1(1/2+h)-F_1(1/2-h)}{2h}.
\end{equation}
Note that if either $f_1$ or $f_l'$ is Lipschitz on $[1/2-\delta_1,1/2+\delta_1]$, then in particular, the function $f_1$ is bounded on this interval by a constant $\kappa>0$, and consequently, the mean value theorem implies 
 \begin{equation*}
|\textsc{cdq}(h)| \leq \kappa,
 \end{equation*} 
 when $h\in (0,\delta_1)$.
Finally, it remains to handle the two integrals on the right side of line~\eqref{threepieces}. We only bound the first one, since the second is essentially the same. Noting that $h\in (0,\delta_1)$, we split the integral into two pieces over $I_1:=[0,1/2-\delta_1]$ and $I_2:=[1/2-\delta_1,1/2-h]$. For any $\theta\in I_1$, we have $\exp(-2t(1/2-h-\theta)^2) \leq \exp(-2t(\delta_1-h)^2)$. Meanwhile, over the second interval, $f_1$ is bounded by a constant $\kappa$. It follows that,
\footnotesize
\begin{equation}\label{lastforfirstlemma}
\ts\frac{1 }{2 h}\displaystyle \int_0^{1/2-h} e^{-2t\,(1/2-h-\theta)^2}f_1(\theta)d\theta \ \leq \  
 \ts\frac{1}{2h}e^{-2t(\delta_1-h)^2}F_1(1/2-\delta_1)+ \ts\frac{\kappa}{2h} \displaystyle\int_{1/2-\delta_1}^{1/2-h} e^{-2t(1/2-h-\theta)^2}d\theta.
 \end{equation}
 \normalsize
The first term on the right is clearly at most $\ts\frac{\kappa}{h\sqrt{t}}$ for some constant $\kappa$. Also, the second term can be calculated exactly as $\ts\frac{1 }{4h\sqrt{t}}\cdot \sqrt{\pi/2}\,\text{Erf}(\sqrt{2}(\delta_1-h)\sqrt{t})$, where $\text{Erf}$ is the error function defined by 
\begin{equation}\label{erfdef}
 \text{Erf}(r):=\ts\frac{2}{\sqrt{\pi}}\int_0^re^{-s^2}ds,
\end{equation}
which satisfies $\text{Erf}(r)\leq 1$ for all real numbers $r$. Hence, the both terms on the right side of line~\eqref{lastforfirstlemma} are most $\kappa/h\sqrt{t}$ for some constant $\kappa$.
\qed

\subsubsection{Proof of Lemma~\ref{crosslemma}}\label{sec:crosslemma}
Note that the quantity $\ts\frac{2}{m_1^2}\binom{m_1}{2}$ in the statement of the lemma is at most 1 for all $m_1$.
Expanding the product $\Delta_1\Delta_2$ gives
\small
\begin{equation}\label{lemma3}
\begin{split}
\E[\Delta_1\Delta_2|\mathcal{D}] &=\ts\frac{1}{h^2}\,\E\Big[K\Big(\ts\frac{\bar{Q}_{t}(\tilde{X}_{1,1})-1/2}{h}\Big)\cdot K\Big(\ts\frac{\bar{Q}_{t}(\tilde{X}_{2,1})-1/2}{h}\Big)\Big|\D\Big] \\[0.2cm]
& \ \ \ \ \ \  \ -2f_1(\ts\frac{1}{2})\cdot \ts\frac{1}{h}\E\Big[K\Big(\ts\frac{\bar{Q}_{t}(\tilde{X}_{1,1})-1/2}{h}\Big)\Big|\D\Big] +f_1(\ts\frac{1}{2})^2.
\end{split}
\end{equation}
\normalsize
We proceed by analyzing the first two terms on the right side separately.
To handle the second term, we will bound it from below, because its negative contribution will be needed to obtain the stated result. Recall the function $\varphi_t$ from line~\eqref{varphiprob} in the proof of Lemma~\ref{squarelemma}, which satisfies
\begin{equation}\label{secondKidentity}
\ts\frac{1}{h}\E\Big[K\Big(\ts\frac{\bar{Q}_{t}(\tilde{X}_{1,1})-1/2}{h}\Big)\Big|\mathcal{D}\Big]=\ts\frac{1}{2h}\displaystyle\int_0^1\varphi_t(\theta;h)f_1(\theta)d\theta.
\end{equation}
Define the following ``lower-bounding function'', $\ell_t:[0,1]\to\R$, which acts an approximate indicator on $[1/2-h,1/2+h]$,
\begin{equation*}
\ell_t(\theta;h):=
\begin{cases}
& 1 - e^{-2t(1/2+h-\theta)^2}-e^{-2t(1/2-h-\theta)^2}, \text{ \ \ for } \theta\in [\ts\frac{1}{2}-h,\ts\frac{1}{2}+h],\\
& 0, \text{  \ \ \ \ \ \ \ \ \ \ \ \ \ \ \   \ \ \ \ \  \ \ \ \ \ \  \ \ \ \ \ \  \ \ \ \ \ \ \ \ \ \ \ \ \ \ \ otherwise.}
\end{cases}
\end{equation*}
\normalsize
Then, Hoeffding's inequality implies
\begin{equation*}
\varphi_t(\theta;h) \geq \ell_t(\theta;h),
\end{equation*}
for all $\theta\in[0,1]$, and all $h\in(0,\delta_1)$. The next step is to integrate this bound over $[0,1]$. In doing this, note that the Lipschitz condition on either $f_1$ or $f_1'$ implies that $f_1$ is bounded by a constant $\kappa$ on $[1/2-h,1/2+h]$. Consequently,
\begin{equation}\label{cdqlower}
\footnotesize
\begin{split}
\ts\frac{1}{h}\E\Big[K\Big(\ts\frac{\bar{Q}_{t}(X_{1,1})-1/2}{h}\Big)\Big|\mathcal{D}\Big] &\geq  \ts\frac{F_1(1/2+h)-F_1(1/2-h)}{2h} \\[0.2cm]
 & \ \ \ \ \ \ \ -\ts\frac{\kappa}{2h}\displaystyle \int_{1/2-h}^{1/2+h}\Big(e^{-2t(1/2+h-\theta)^2}+e^{-2t(1/2-h-\theta)^2}\Big)d\theta\\[0.3cm]
&=  \ \textsc{cdq}(h)- \kappa\cdot \ts\frac{\sqrt{\pi/2}}{2h\sqrt{t}} \cdot \text{Erf}(2\sqrt{2}h\sqrt{t})\\[0.3cm]
&\geq \textsc{cdq}(h)-\ts\frac{\kappa}{h\sqrt{t}},
\end{split}
\end{equation}
\normalsize
where  $\textsc{cdq}(h)$ was defined in line~\eqref{cdqdef}, and the second line follows from an exact calculation using the error function $\text{Erf}$ defined in line~\eqref{erfdef}, which is always bounded in magnitude by 1.

To handle the first term in line~\eqref{lemma3}, note that for the rectangular kernel, we have the inequality $K(a)K(b)\leq \ts\frac{1}{2}\min\{K(a),K(b)\}$ for all $a,b\in\R$. So, using the identity~\eqref{secondKidentity} and the bound~\eqref{uupper} gives
\scriptsize
\begin{equation}\label{doubleint}
\begin{split}
\ts\frac{1}{h^2}\,\E\Big[K\Big(\ts\frac{\bar{Q}_{t}(\tilde{X}_{1,1})-1/2}{h}\Big)\cdot K\Big(\ts\frac{\bar{Q}_{t}(\tilde{X}_{2,1})-1/2}{h}\Big)\Big|\D\Big] 
&=\int_{\mathcal{X}} \int_{\mathcal{X}} \ts\frac{1}{h^2}\,\E\bigg[ K\Big(\ts\frac{\bar{Q}_{t}(x)-1/2}{h}\Big) K\Big(\ts\frac{\bar{Q}_{t}(x')-1/2}{h}\Big)\Big|\D\bigg] d\mu_1(x)d\mu_1(x') \\[0.3cm]
&\leq \int_{\mathcal{X}} \int_{\mathcal{X}} \ts\frac{1}{4h^2}\cdot \min\Big\{ \varphi_t(\theta(x);h) \, , \, \varphi_t(\theta(x');h)\Big\} d\mu_1(x)d\mu_1(x') \\[0.3cm]
&=\int_0^1 \int_0^1 \ts\frac{1}{4h^2}\cdot \min\Big\{ \varphi_t(\theta;h) \, , \, \varphi_t(\theta';h)\Big\} f_1(\theta) f_1(\theta')\, d\theta\, d\theta'\\[0.3cm]
&\leq \int_0^1 \int_0^1 \ts\frac{1}{4h^2}\cdot \min\Big\{ \mathscr{U}_t(\theta;h) \, , \, \mathscr{U}_t(\theta';h) \Big\} f_1(\theta) f_1(\theta')\, d\theta\, d\theta'.
\end{split}
\end{equation}
\normalsize
Let the last integral on the right side be denoted by $(**)$. We will bound this integral by decomposing the unit square $[0,1]^2$ into the five regions displayed below in Figure~\ref{fig:square}, where $\theta$ is measured along the x-axis and $\theta'$ is measured along the y-axis. The decomposition gives
\small
\begin{align*}
(**) \leq \int_{1/2-h}^{1/2+h}\int_{1/2-h}^{1/2+h} \ts\frac{1}{4h^2}f_1(\theta)f_1(\theta')d\theta d\theta'  & 
+ \int_{\text{A}} \ts\frac{1 }{4h^2}e^{-2t\,(1/2-h-\theta)^2}f_1(\theta) f_1(\theta') d\theta d\theta'\\[0.2cm]
& + \int_{\text{B}} \ts\frac{1 }{4h^2}e^{-2t\,(1/2-h-\theta')^2} f_1(\theta) f_1(\theta')d\theta d\theta'\\[0.2cm]
& + \int_{\text{C}} \ts\frac{1}{4h^2}e^{-2t\,(\theta-(1/2+h))^2}f_1(\theta) f_1(\theta') d\theta d\theta'\\[0.2cm]
& + \int_{\text{D}} \ts\frac{1}{4h^2}e^{-2t\,(\theta'-(1/2+h))^2}f_1(\theta) f_1(\theta') d\theta d\theta'.
\end{align*}
\normalsize
\vspace{-0.5cm}
\begin{figure}[h!]
\centering
{\includegraphics[angle=-90,
  width=0.5\linewidth]{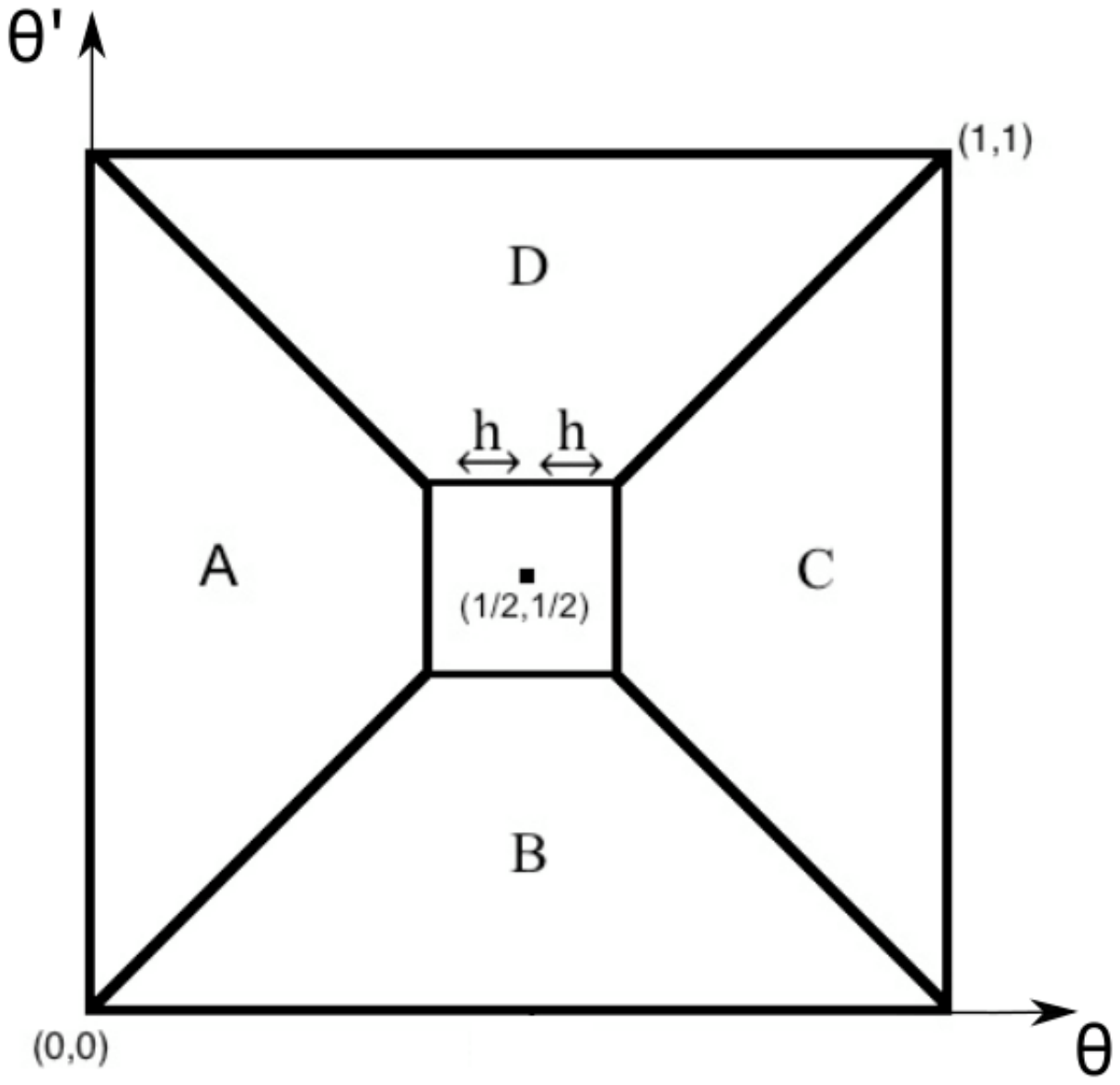}}
\caption{The unit square $[0,1]^2$ is decomposed into five regions.
}
\label{fig:square}
\end{figure}

\noindent By symmetry, each of the integrals over A, B, C, and D are equal, and so we only bound the one over A. At the end of this subsection, we will give the details for deriving the following bound,
\begin{equation}\label{intA}
\begin{split}
\int_{\text{A}} \ts\frac{1}{4h^2}e^{-2t\,(1/2-h-\theta)^2} f_1(\theta) f_1(\theta') d\theta d\theta' 
&\leq \kappa\Big(\ts\frac{1}{h\sqrt{t}}+\ts\frac{1}{h^2t}\Big),
\end{split}
\end{equation}
for some constant $\kappa>0$ not depending on $t$ or $m_1$. Lastly, it is simple to check that the integral over the square region $[1/2-h,1/2+h]^2$
is given by
\begin{equation}\label{cdqint}
\begin{split}
\int_{1/2-h}^{1/2+h}\int_{1/2-h}^{1/2+h} \ts\frac{1}{4h^2}f_1(\theta)f_1(\theta')d\theta d\theta' = \textsc{cdq}(h)^2.
\end{split}
\end{equation}
 Combining the work from lines~\eqref{lemma3},~\eqref{cdqlower},~\eqref{doubleint},~\eqref{intA}, and~\eqref{cdqint} gives
 \begin{equation}\label{dqsquared}
 \begin{split}
 \ts\frac{1}{m_1}\binom{m_1}{2} \E[\Delta_1\Delta_2|\D] &\leq 
 \Big(\textsc{cdq}(h)^2+\ts\frac{\kappa}{h\sqrt{t}}+\ts\frac{1}{h^2t}\Big) -2f_1(\ts\frac{1}{2})\cdot \Big( \textsc{cdq}(h)-\ts\frac{\kappa}{h\sqrt{t}}\Big)+f_1(\ts\frac{1}{2})^2\\[0.2cm]
 &= \Big(\textsc{cdq}(h)-f_1(\ts\frac{1}{2})\Big)^2+\kappa\Big(\ts\frac{1}{h\sqrt{t}}+\ts\frac{1}{h^2t}\Big),
 \end{split}
 \end{equation}
 where we have absorbed a factor of $f_1(\ts\frac{1}{2})$ into $\kappa$.
Regarding the difference quotient, there is a constant $\kappa>0$, such that whenever $h\in (0,\delta_1)$,
 \begin{align}\label{splitbounds}
 |\textsc{cdq}(h)-f_1(\ts\frac{1}{2})|\leq 
 \begin{cases}
 & \kappa h, \text{ \ \ \  when $f_1$ is Lipschitz on $[1/2-\delta_1,1/2+\delta_1]$,}\\
 & \kappa h^2 \, \text{ \ \  when $f_1'$ is Lipschitz on $[1/2-\delta_1,1/2+\delta_1]$.}
 \end{cases}
 \end{align}
 Inserting these bounds into line~\eqref{dqsquared} yields the statement of the lemma. Below, we give two paragraphs providing the details for the previous line~\eqref{splitbounds} and the bound~\eqref{intA}.\qed
 ~\\

 \paragraph{Details for line~\eqref{splitbounds}}

 To give some detail for the bounds in line~\eqref{splitbounds}, first note that
 \begin{equation}\label{cdqintegral}
 \textsc{cdq}(h)-f_1(\ts\frac{1}{2}) =\ts\frac{1}{2h} \displaystyle\int_{1/2-h}^{1/2+h} (f_1(\theta)-f_1(\ts\frac{1}{2}))d\theta.
 \end{equation}
 When $f_1$ is Lipschitz on $[1/2-\delta_1,1/2+\delta_1]$, the integrand above satisfies $|f_l(\theta)-f_1(\ts\frac{1}{2})|\leq \kappa h$, which leads to the first bound in line~\eqref{splitbounds}. To handle the case when $f_l'$ is Lipschitz, consider the Taylor expansion
\begin{equation}\label{expansion}
f_1(\theta)-f_1(\ts\frac{1}{2}) = (\theta-\ts\frac{1}{2})f_1'(\ts\frac{1}{2})+R(\theta),
\end{equation}
where, for $\theta\in [1/2-h,1/2+h]$, the remainder is defined by
\begin{equation*}
R(\theta):=\int_0^1 \Big[f_1'(\ts\frac{1}{2}+s(\theta-\ts\frac{1}{2}))-f_1'(\ts\frac{1}{2})\Big](\theta-\ts\frac{1}{2})ds.
\end{equation*}
When $f_l'$ is Lipschitz on $[1/2-\delta_1,1/2+\delta_1]$, there is a constant $\kappa>0$ such  that whenever $\theta\in[1/2-h,1/2+h]$, the remainder will satisfy
\begin{equation*}
\begin{split}
|R(\theta)| & \leq (\theta-\ts\frac{1}{2})^2 \int_0^1  sds\\[0.2cm]
&\leq \kappa h^2.
\end{split}
\end{equation*}
 Combining the previous step with lines~\eqref{expansion} and~\eqref{cdqintegral} proves the second case in line~\eqref{splitbounds}, since the term $(\theta-\ts\frac{1}{2})f_1'(\ts\frac{1}{2})$ vanishes under the integral in line~\eqref{cdqintegral}.

\paragraph{Details for the bound in line~\eqref{intA}} 
Define the quantity
$$(***):=\int_{\text{A}} \ts\frac{1}{4h^2}e^{-2t\,(1/2-h-\theta)^2} f_1(\theta) f_1(\theta') d\theta d\theta' .$$
By direct calculation, 
\begin{equation*}
\begin{split}
(***) &
 =\int_0^{1/2-h}  \ts\frac{1}{4h^2}e^{-2t\,(1/2-h-\theta)^2} \cdot f_1(\theta) \cdot \displaystyle\bigg(\int_{\theta}^{1-\theta}  f_1(\theta')d\theta'\bigg)d\theta\\[0.3cm]
&=\int_0^{1/2-h}  \ts\frac{1}{4h^2}e^{-2t\,(1/2-h-\theta)^2} f_1(\theta) \cdot \big(F_1(1-\theta)-F_1(\theta)\big)d\theta.
\end{split}
\end{equation*}
\normalsize
Next, because $h<\delta_1$, we may split the interval $[0,1/2-h]$ into the union of $I_1:=[0,1/2-\delta_1]$ and $I_2:=[1/2-\delta_1, 1/2-h]$. When $f_1$ or $f_1'$ is Lipschitz on $[1/2-\delta_1,1/2+\delta_1]$, it follows in particular that $f_1$ is bounded by a positive constant $\kappa$ on $I_2$. Similarly, the mean value theorem implies that $|F_1(1-\theta)-F_1(\theta)|\leq \kappa (1-2\theta)$ on $I_2$.

 Alternatively, on the interval $I_1$, we note that $e^{-2t(1/2-h-\theta)^2}$ is at most $e^{-2t(\delta_1-h)^2}$, and also the function $|F_1(1-\theta)-F_1(\theta)|$ is clearly at most 1. Combining these observations, it follows that
\begin{align}\label{secondtolast}
(***) & 
 \ \leq \  \ts\frac{1}{4h^2}e^{-2t\,(\delta_1-h)^2}F_1(1/2-\delta_1) + \ts\frac{\kappa}{4h^2}\displaystyle \int_{1/2-\delta_1}^{1/2-h} (1-2\theta)e^{-2t(1/2-h-\theta)^2} \, d\theta.
\end{align}
\normalsize
Finally, the last integral on the right side can be calculated exactly with the change of variable $v:=1/2-h-\theta$, which gives
\begin{equation}\label{lastone}
\begin{split}
\ts\frac{\kappa}{4h^2}\displaystyle \int_{1/2-\delta_1}^{1/2-h} (1-2\theta)e^{-2t(1/2-h-\theta)^2} \, d\theta  \ & = \ \kappa \int_{0}^{\delta_1-h} \ts\frac{2h+2v}{4h^2} e^{-2tv^2}dv\\[0.3cm]
& = \frac{1-e^{-2t(\delta_1-h)^2}}{8h^2t}+\frac{\sqrt{2\pi}\, \text{Erf}(\sqrt{2t}(\delta_1-h))}{8 h\sqrt{t}},
\end{split}
\end{equation}
where the error function $\text{Erf}$ is defined in line~\eqref{erfdef}, and is always bounded in magnitude by 1. Combining lines~\eqref{secondtolast} and~\eqref{lastone} shows that up to a positive constant $\kappa$ not depending on $t$ or $m_1$, we have
$$(***) \leq \kappa\Big(\ts\frac{1}{ht^2}+\ts\frac{1}{h\sqrt{t}}\Big),$$
as needed.

\normalsize

\section{Validation of assumption \textbf{A2}}

To validate assumption \textbf{A2}, we carried out the following experiments to see how well the distribution $\mathcal{L}(\theta(X)|\D,Y=l)$ can be approximated by a distribution that is known to satisfy~\textbf{A2}. (Note that a similar set of validation experiments has also been presented in the supplement of~\cite{lopes2019} with different datasets.) Since the random variable $\theta(X)$ takes values in the interval [0,1], a familiar class of distributions to consider is the Beta family. Recall that for any $\tau\in[0,1]$, the density $g(\tau;\alpha,\beta)$ of the Beta($ \alpha,\beta$) distribution is given by
$$g(\tau; \alpha,\beta) := \frac{1}{B(\alpha,\beta)} \tau^{\alpha-1}(1-\tau)^{1-\beta},$$
where $\alpha,\beta>0$ are parameters, and $B(\alpha,\beta)$ is the Beta function.
From this formula, it is clear that for any choice the parameters, the function $g(\cdot;\alpha,\beta)$ is Lipschitz in a neighborhood of 1/2.

The first step in these experiments involved generating approximate samples from the distribution $\mathcal{L}(\theta(X)|\D,Y=l)$. For this purpose, we prepared the training set $\mathcal{D}$ and ground truth set $\mathcal{D}_{\text{ground}}$ for both the `abalone' and `landsat satellite' data, as described in Section 4 of the main text. To approximate the function $\theta$, we trained a very large ensemble of 10,000 decision trees $\{Q_i\}$ via random forests (abbrev. `RF') on $\D$, and then used the approximation $\theta(\cdot)\approx \bar{Q}_t(\cdot)$. Letting $X_{1,l}',\dots,X_{n'_l,l}'$ denote the samples in $\mathcal{D}_{\text{ground}}$ from class $l$,  we used the values $\bar{Q}_t(X_{1,l}'),\dots,\bar{Q}_t(X_{n'_l,l}')$ as approximate samples from  $\mathcal{L}(\theta(X)|\D,Y=l)$. (Note that $n_0'+n_1'=|\mathcal{D}_{\text{ground}}|$.)

In the case of the `landsat satellite' data, $|\mathcal{D}_{\text{ground}}|=3,217$, and in the case of the `abalone' data, $|\mathcal{D}_{\text{ground}}=2,088$. Hence, for both datasets, a fairly large number of approximate samples from $\mathcal{L}(\theta(X)|\D,Y=l)$ were available. These approximate samples were then used in conjunction with the method of moments to estimate the parameters $\alpha$ and $\beta$. The implementation of the method of moments was done with the ``mme'' option in the R package `fitdistrplus', and below, we denote the estimates associated with class $l$ as $\hat{\alpha}_l$ and $\hat{\beta}_l$.

To see how well the Beta distributions fit the data, we used quantile-quantile (QQ) plots as a diagonostic. Specifically, we sorted the values $\bar{Q}_t(X_{1,l}'),\dots,\bar{Q}_t(X_{n_l,l}')$ and plotted them against a corresponding set of quantiles from the Beta$(\hat{\alpha}_l,\hat{\beta}_l)$ distribution. The results are displayed in purple in the figures below, and we can see that there is a good overall conformity with the diagonal line.
 We have also marked a horizontal line at 1/2 to clarify that the fit is good in a neighborhood of 1/2 (i.e.~the region of interest for \textbf{A2}).

Next, we turn our attention to the small deviations of the purple curves from the diagonal line, which occur mostly in the tails of the distributions. To look at this more carefully, we generated an i.i.d.~sample of size $n_l'$ from the fitted distribution Beta$(\hat{\alpha}_l,\hat{\beta}_l)$, and plotted the sorted samples against the same quantiles used previously  (but with the results shown in green). 

The importance of the green curves is that they show that some deviation can be expected from the diagonal line, even when the samples and quantiles are both obtained from the same distribution (i.e. Beta$(\hat{\alpha}_l,\hat{\beta}_l)$). Moreover, the green curves have deviations from the diagonal that resemble the deviations of the purple curves. Consequently, the green curves give further support to the conclusion that $\mathcal{L}(\theta(X)|\mathcal{D},Y=l)$ is well approximated by a Beta distribution.

\begin{figure*}[h!]
\centering
\centering
  {\includegraphics[angle=0,
  width=0.4\linewidth]{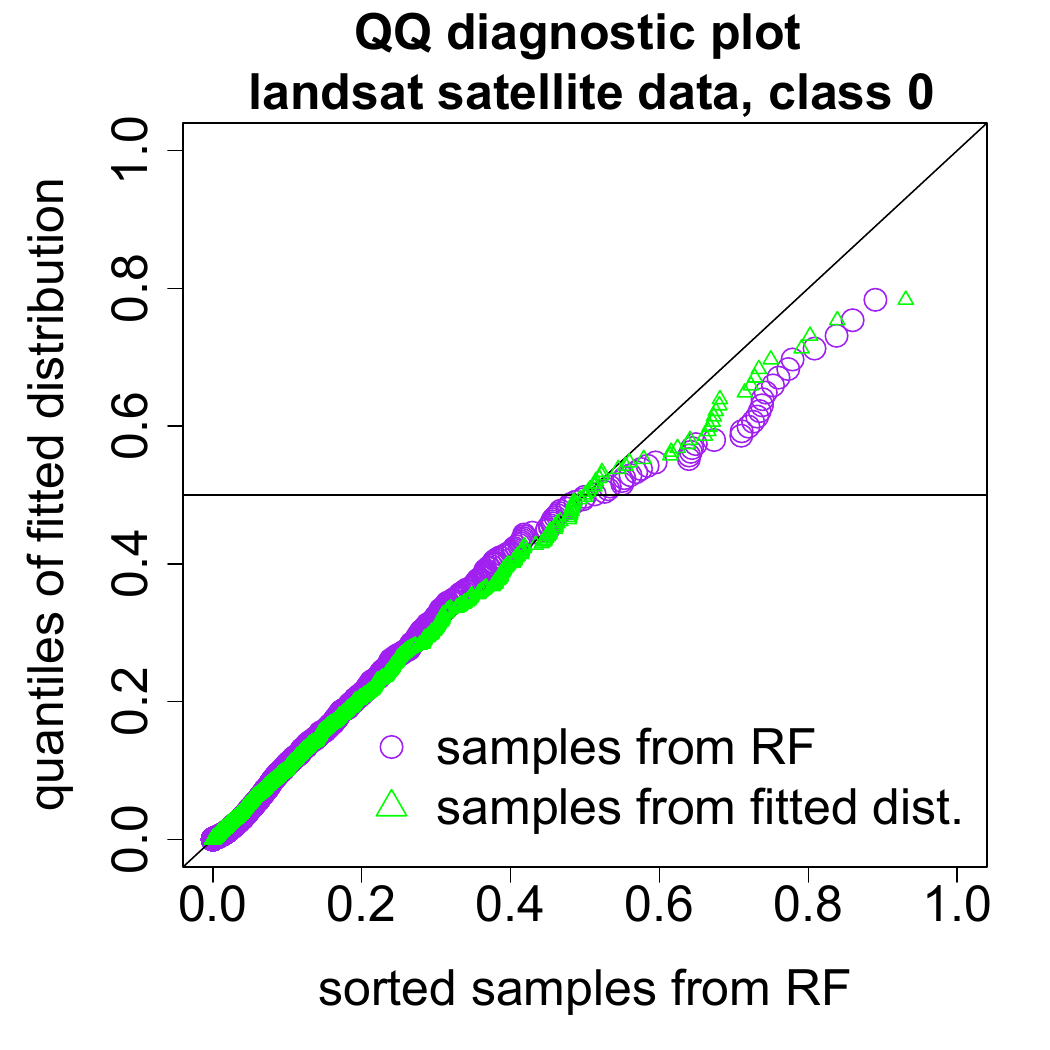}}
  \ \ \ \ 
  {\includegraphics[angle=0,
  width=0.4\linewidth]{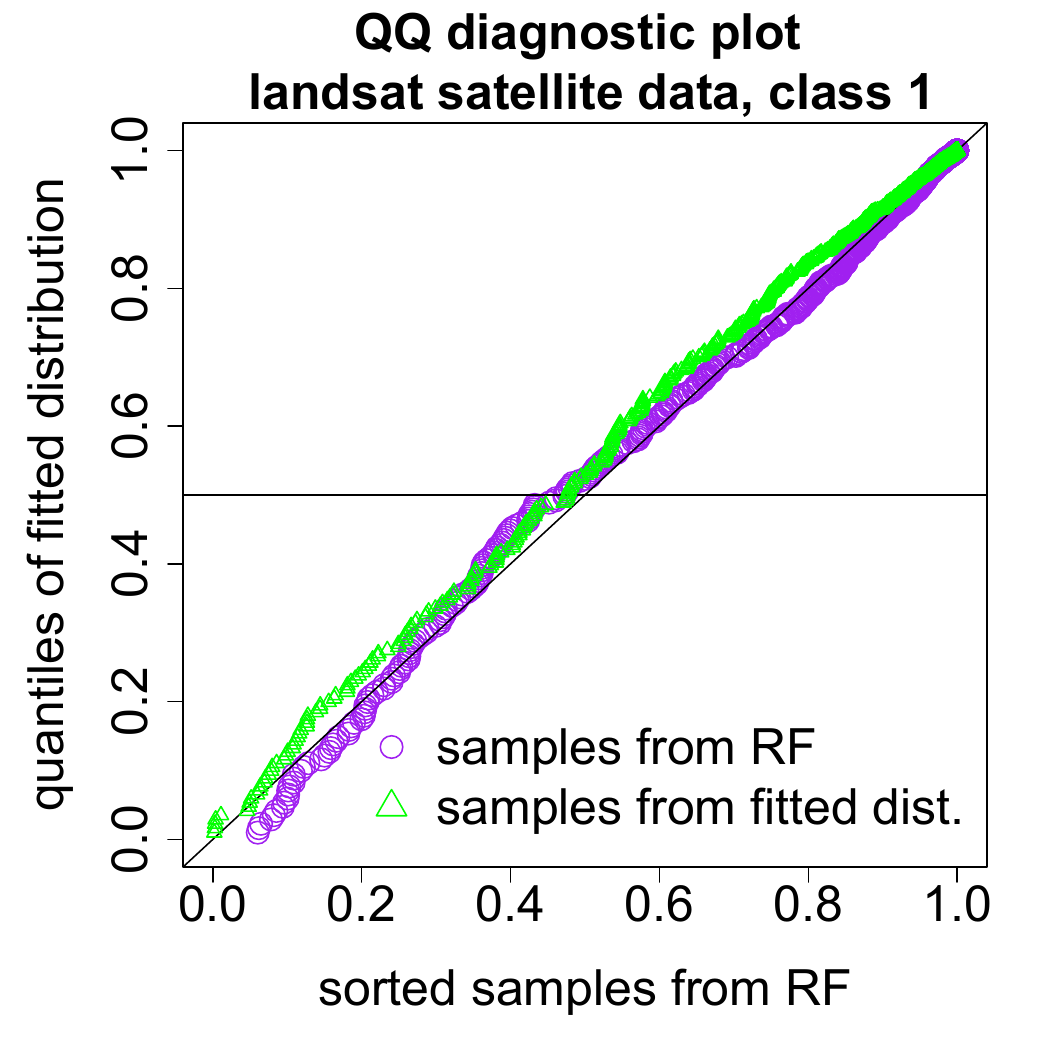}}
  ~\\
  ~\\
  ~\\
  ~\\
  {\includegraphics[angle=0,
  width=0.4\linewidth]{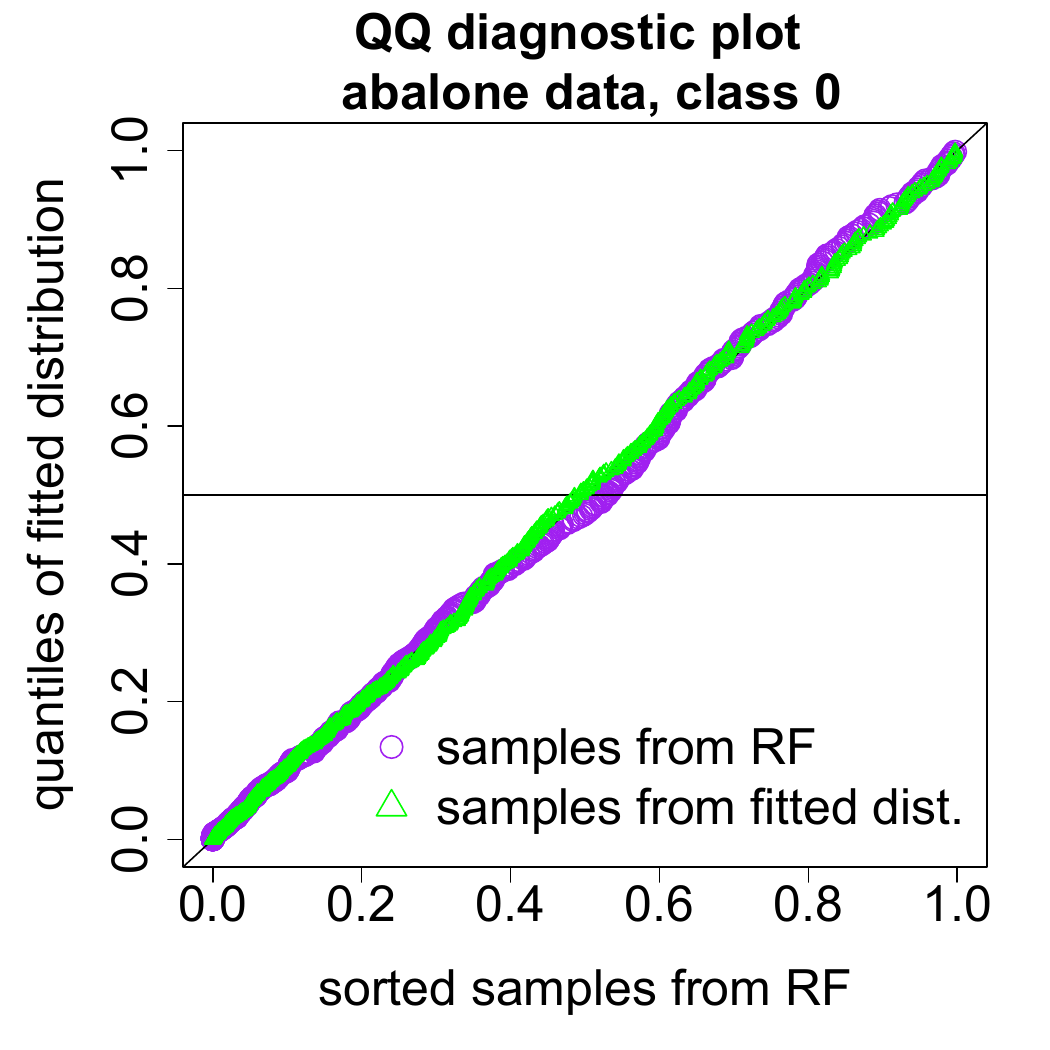}}
  \ \ \ \ 
  {\includegraphics[angle=0,
  width=0.4\linewidth]{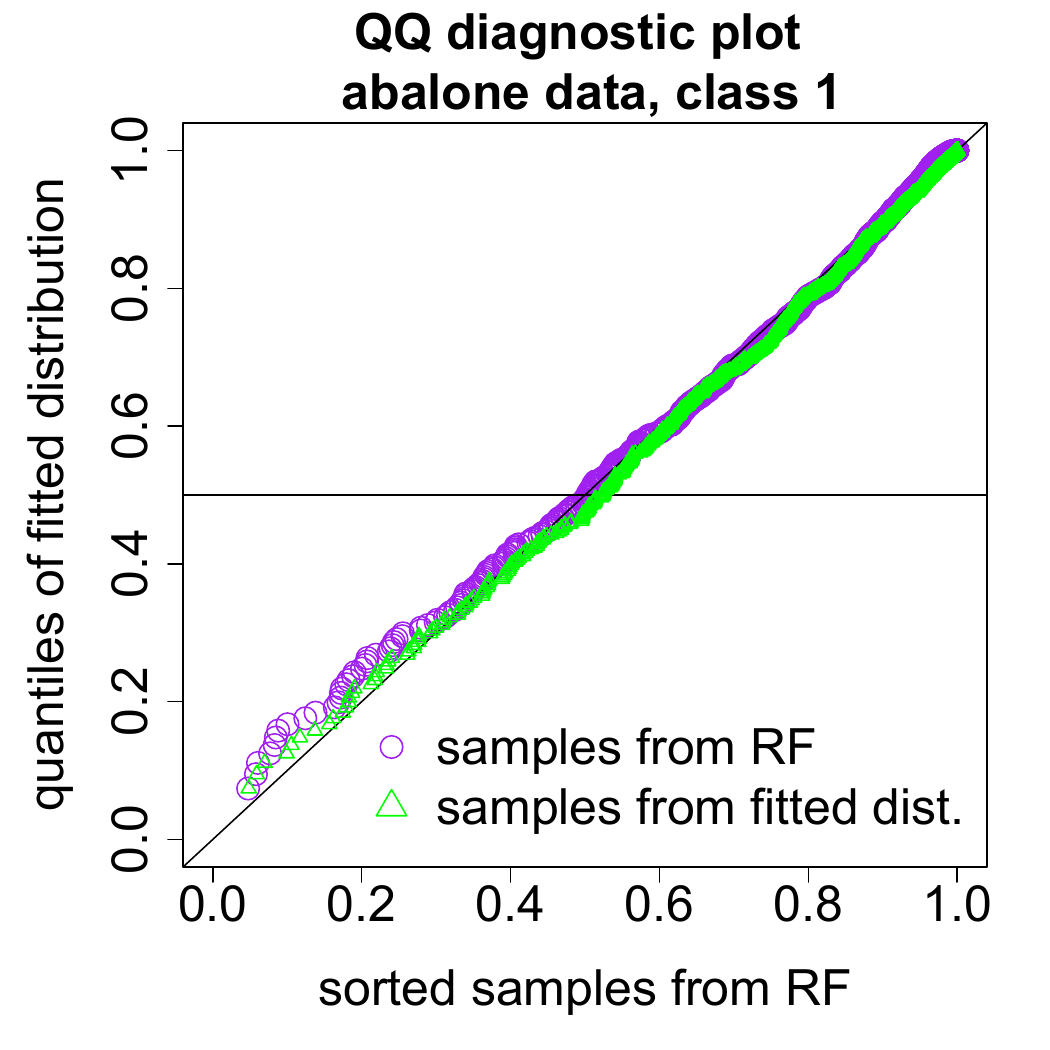}}
\caption{QQ plots for validating assumption \textbf{A2} in the case of the `landsat satellite' and `abalone' data. Overall, the purple curves align well with the diagonal line, indicating a good fit. Note that even when the purple curves deviate from the diagonal line, the deviations are roughly similar to those of a sample drawn from the corresponding fitted distributions (as shown in green). Hence, the small deviations of the purple curves do not necessarily reflect a poor fit.}
\end{figure*}

\section{Comments on data preparation}\label{app:data}

When  subsets of a common dataset were partitioned in separate files in the UCI repository, we combined these subsets before further processing. To make the datasets compatible with binary classification, some label classes were pooled in a few cases. For `abalone', labels 1-8 were set to 0 with others set to 1, for `digits', labels 0-4 were set to 0 with others set to 1, and for `landsat satellite', labels 1-3 were set to 0 with others set to 1. Also, in the `occupancy detection' dataset, 10\% of the labels were randomly selected and then flipped, so that the classification problem would lead to non-trivial error rates.




\begin{thebibliography}{10}

\bibitem{cannings2017}
T.~I. Cannings and R.~J. Samworth.
\newblock Random projection ensemble classification (with discussion).
\newblock {\em Journal of the Royal Statistical Society Series B}, 2017.

\bibitem{lopes2019}
M.~E. Lopes.
\newblock Estimating the algorithmic variance of randomized ensembles via the
  bootstrap.
\newblock {\em The Annals of Statistics}, 47(2):1088--1112, 04 2019.

\bibitem{breiman1996}
L.~Breiman.
\newblock Bagging predictors.
\newblock {\em Machine learning}, 24(2):123--140, 1996.

\bibitem{breiman2001}
L.~Breiman.
\newblock Random forests.
\newblock {\em Machine learning}, 45(1):5--32, 2001.

\bibitem{buhlmannyu}
P.~B{\"u}hlmann and B.~Yu.
\newblock Analyzing bagging.
\newblock {\em The Annals of Statistics}, pages 927--961, 2002.

\bibitem{breimanconsistency}
L.~Breiman.
\newblock Consistency for a simple model of random forests.
\newblock {\em Technical report}, 2004.

\bibitem{hallsamworth}
P.~Hall and R.~J. Samworth.
\newblock Properties of bagged nearest neighbour classifiers.
\newblock {\em Journal of the Royal Statistical Society: Series B},
  67(3):363--379, 2005.

\bibitem{linjeon}
Y.~Lin and Y.~Jeon.
\newblock Random forests and adaptive nearest neighbors.
\newblock {\em Journal of the American Statistical Association},
  101(474):578--590, 2006.

\bibitem{biau2008}
G.~Biau, L.~Devroye, and G.~Lugosi.
\newblock Consistency of random forests and other averaging classifiers.
\newblock {\em Journal of Machine Learning Research}, 9:2015--2033, 2008.

\bibitem{biau2012}
G.~Biau.
\newblock Analysis of a random forests model.
\newblock {\em Journal of Machine Learning Research}, 98888:1063--1095, 2012.

\bibitem{scornetconsistency}
E.~Scornet, G.~Biau, and J.-P. Vert.
\newblock Consistency of random forests.
\newblock {\em The Annals of Statistics}, 43(4):1716--1741, 2015.

\bibitem{NgJordan}
A.~Y. Ng and M.~I. Jordan.
\newblock Convergence rates of the {V}oting {G}ibbs classifier, with
  application to {B}ayesian feature selection.
\newblock In {\em International Conference on Machine Learning}, pages
  377--384, 2001.

\bibitem{hernandez2013}
D.~Hern{\'a}ndez-Lobato, G.~Mart{\'\i}nez-Mu{\~n}oz, and A.~Su{\'a}rez.
\newblock How large should ensembles of classifiers be?
\newblock {\em Pattern Recognition}, 46(5):1323--1336, 2013.

\bibitem{elements}
J.~Friedman, T.~Hastie, and R.~Tibshirani.
\newblock {\em The {E}lements of {S}tatistical {L}earning}.
\newblock Springer, 2001.

\bibitem{tinkamho1998}
T.~K. Ho.
\newblock The random subspace method for constructing decision forests.
\newblock {\em IEEE transactions on pattern analysis and machine intelligence},
  20(8):832--844, 1998.

\bibitem{dietterich2000}
T.~G. Dietterich.
\newblock An experimental comparison of three methods for constructing
  ensembles of decision trees: Bagging, boosting, and randomization.
\newblock {\em Machine learning}, 40(2):139--157, 2000.

\bibitem{schapire2012boosting}
R.~Schapire and Y.~Freund.
\newblock {\em Boosting: Foundations and Algorithms}.
\newblock The MIT Press, 2012.

\bibitem{latinne}
P.~Latinne, O.~Debeir, and C.~Decaestecker.
\newblock Limiting the number of trees in random forests.
\newblock In {\em Multiple Classifier Systems}, pages 178--187. Springer, 2001.

\bibitem{comet}
J.~D. Basilico, M.~A. Munson, T.~G. Kolda, K.~R. Dixon, and W.~P. Kegelmeyer.
\newblock Comet: A recipe for learning and using large ensembles on massive
  data.
\newblock In {\em IEEE International Conference on Data Mining (ICDM)}, pages
  41--50. IEEE, 2011.

\bibitem{swiss}
A.~G. Schwing, C.~Zach, Y.~Zheng, and M.~Pollefeys.
\newblock Adaptive {R}andom {F}orest" -- how many experts€ to ask before
  making a decision?
\newblock In {\em Computer Vision and Pattern Recognition (CVPR), 2011 IEEE
  Conference on}, pages 1377--1384. IEEE, 2011.

\bibitem{oshiro}
T.~M. Oshiro, P.~S. Perez, and J.~A. Baranauskas.
\newblock How many trees in a random forest?
\newblock In {\em Machine Learning and Data Mining in Pattern Recognition},
  pages 154--168. Springer, 2012.

\bibitem{lopes2013}
M.~E. Lopes.
\newblock The convergence rate of majority vote under exchangeability.
\newblock {\em arXiv:1303.0727}, 2013.

\bibitem{cannings2015}
T.~I. Cannings and R.~J. Samworth.
\newblock Random projection ensemble classification.
\newblock {\em arXiv:1504.04595}, 2015.

\bibitem{lopes2016}
M.~E. Lopes.
\newblock A sharp bound on the computation-accuracy tradeoff for majority
  voting ensembles.
\newblock {\em arXiv:1303.0727}, 2016.

\bibitem{simon}
L.~Simon.
\newblock {\em Lectures on {G}eometric {M}easure {T}heory}.
\newblock The Australian National University, Mathematical Sciences Institute,
  Centre for Mathematics \& its Applications, 1983.

\bibitem{cucker}
P.~B{\"u}rgisser and F.~Cucker.
\newblock {\em Condition: The geometry of numerical algorithms}, volume 349.
\newblock Springer, 2013.

\bibitem{bujastuetzle2006}
A.~Buja and W.~Stuetzle.
\newblock Observations on bagging.
\newblock {\em Statistica Sinica}, pages 323--351, 2006.

\bibitem{efron2014}
B.~Efron.
\newblock Estimation and accuracy after model selection.
\newblock {\em Journal of the American Statistical Association},
  109(507):991--1007, 2014.

\bibitem{slepian1962}
D.~Slepian.
\newblock The one-sided barrier problem for {G}aussian noise.
\newblock {\em Bell System Technical Journal}, 41(2):463--501, 1962.

\bibitem{sidak1968}
Z.~Sid{\'a}k.
\newblock On multivariate normal probabilities of rectangles: their dependence
  on correlations.
\newblock {\em The Annals of Mathematical Statistics}, pages 1425--1434, 1968.

\bibitem{schervish2012}
M.~J. Schervish.
\newblock {\em Theory of Statistics}.
\newblock Springer, 2012.

\bibitem{tsybakov}
A.~Tsybakov.
\newblock {\em Introduction to Nonparametric Estimation}.
\newblock Springer, 2009.

\bibitem{meister}
A.~Meister.
\newblock {\em Deconvolution Problems in Nonparametric Statistics}.
\newblock Springer, 2009.

\bibitem{delaigle}
A~Delaigle and A.~Meister.
\newblock Density estimation with heteroscedastic error.
\newblock {\em Bernoulli}, pages 562--579, 2008.

\bibitem{hesse}
C.~Hesse.
\newblock How many ``good'' observations do you need for" fast" density
  deconvolution from supersmooth errors.
\newblock {\em Sankhy{\=a}: The Indian Journal of Statistics, Series A}, pages
  491--506, 1996.

\bibitem{JMVAdecon}
X.-F. Wang and D.~Ye.
\newblock Conditional density estimation in measurement error problems.
\newblock {\em Journal of Multivariate Analysis}, 133:38--50, 2015.

\bibitem{meister2010}
A.~Meister and M.~H. Neumann.
\newblock Deconvolution from non-standard error densities under replicated
  measurements.
\newblock {\em Statistica Sinica}, pages 1609--1636, 2010.

\bibitem{uci}
M.~Lichman.
\newblock {UCI} machine learning repository, 2013.

\bibitem{randomForestsCitation}
A.~Liaw and M.~Wiener.
\newblock Classification and regression by randomforest.
\newblock {\em R News}, 2(3):18--22, 2002.

\bibitem{vaart}
A.~W. van~der Vaart.
\newblock {\em Asymptotic Statistics}.
\newblock Cambridge University Press, 2000.

\bibitem{kallenberg}
O.~Kallenberg.
\newblock {\em Foundations of Modern Probability}.
\newblock Springer, 2006.

\bibitem{boucheron2013}
S.~Boucheron, G.~Lugosi, and P.~Massart.
\newblock {\em Concentration Inequalities: A Nonasymptotic Theory of
  Independence}.
\newblock Oxford University Press, 2013.

\end{thebibliography}


\end{document}